\documentclass[12pt]{article}
\usepackage{amssymb,a4}
\usepackage{amsmath,amsfonts,amssymb,amsthm}
\newcommand{\f}{\mathfrak}

\newcommand{\al}{\alpha}

\def\Res{{\rm Res}}
\def\wt{{\rm wt}}
\def\De{\Delta}
\def\de{\delta}
\def\be{\beta}

\newcommand{\la}{\lambda}\newcommand{\La}{\Lambda}

\def\C{{\mathbb C}}

\def\Z{{\mathbb Z}}

\def\1{{\bf 1}}
\def \End{{\rm End}}
\def \Hom{{\rm Hom}}
\def \Vir{{\rm Vir}}
\def \pf{\noindent {\bf Proof: \,}}
\def\theequation{5.\arabic{equation}}
\def \a{\alpha}
\def \b{\beta}
\def \e{\epsilon}
\def \h{\mathfrak{h}}
\def \l{\lambda}
\def \w{\omega}

\renewcommand{\theequation}{\thesection.\arabic{equation}}

\begin{document}
\newcommand{\Free}[1]{{M(1)}^{#1}}
\newcommand{\Fremo}[1]{M(1,#1)}
\newcommand{\charge}[1]{V_{L}^{#1}}
\newcommand{\charhalf}[2]{V_{#1+L}^{#2}}
\newcommand{\charlam}[1]{V_{#1+L}}
\newcommand{\del}{\partial}

\newcommand{\namelistlabel}[1]{\mbox{#1}\hfill}
\newtheorem{theorem}{Theorem}[section]
\newtheorem{prop}[theorem]{Proposition}
\newtheorem{lem}[theorem]{Lemma}
\newtheorem{coro}[theorem]{Corollary}
\newtheorem{remark}[theorem]{Remark}
\theoremstyle{definition}
\newtheorem{definition}[theorem]{Definition}
\newtheorem{example}[theorem]{Example}
\newenvironment{namelist}[1]{%
\begin{list}{}
{
\let\makelabel\namelistlabel
\settowidth{\labelwidth}{#1}
 \setlength{\leftmargin}{1.1\labelwidth}
 }
 }{%
 \end{list}}
 \begin{center}
{\Large {\bf Rationality of vertex operator algebras}} \\

\vspace{0.5cm} Chongying Dong\footnote{Supported by NSF grants,
China NSF grant 10328102 and  a Faculty research grant from  the
University of California at Santa Cruz.}
\\
Department of Mathematics\\ University of
California\\ Santa Cruz, CA 95064 \\
Cuipo Jiang\footnote{Supported  by China NSF grant
10571119.}\\
 Department of Mathematics\\ Shanghai Jiaotong University\\
Shanghai 200030 China\\
\end{center}
\hspace{1cm}
\begin{abstract} It is shown that a simple vertex operator algebra
$V$ is rational if and only if its Zhu algebra $A(V)$ is
semisimple and each irreducible admissible $V$-module is ordinary.
A contravariant form on a Verma type admissible $V$-module is
constructed and the radical is exactly the maximal proper
submodule. As an application the rationality of $V_L^+$ for any
positive definite even lattice is obtained.
\end{abstract}
 \section{Introduction}
\def\theequation{1. \arabic{equation}}
\setcounter{equation}{0}

One of the most important problems in the representation theory of
vertex operator algebras is to determine various module
categories. Although there are several  ways to define modules for
a vertex operator algebra for different purposes, there are
essentially three different notions of modules, namely,
$${\rm weak\ modules \supset admissible\ modules \supset ordinary\ modules.}$$
The ordinary module was first defined in \cite{FLM} in the
construction of moonshine vertex operator algebra $V^{\natural}.$
The ordinary modules are graded by the eigenvalues of the degree
operator $L(0)$ with finite dimensional eigenspaces and the
eigenvalues are bounded below. In order to study the modular
invariance of trace functions for vertex operator algebras, the
admissible module was introduced in \cite{Z} ( see also
\cite{DLM2}). The admissible modules are $\Z_+$-graded but the
homogenous spaces are not necessarily finite dimensional and
$L(0)$ is not assumed to be semisimple.  The weak modules do not
have any grading restriction.

There are two basics in understanding module categories. The first
one is the classification of irreducible objects and the other is
 whether or not the category is semisimple. For vertex operator
algebras, the admissible module category is the most important one
for various considerations. We call a vertex operator algebra {\em
rational} if the admissible module category is semisimple
\cite{DLM2}. This rationality definition is essentially the
rationality first defined in \cite{Z} with two more assumptions:
(a) There are only finitely many irreducible admissible modules,
(b) each irreducible admissible module is ordinary. It was proved
in \cite{DLM2} that (a) and (b) follow from the semisimplicity of
the admissible module category. Many well known vertex operator
algebras (such as lattice vertex operator algebras \cite{B},
\cite{FLM}, \cite{D1}, \cite{DLM1}, the integrable affine vertex
operator algebras \cite{FZ}, the Virasoro vertex operator algebras
associated to the discrete series \cite{W1}, \cite{DMZ}, the
framed and code vertex operator algebras \cite{D2}, \cite{DGH},
\cite{M}, certain $W$-algebras \cite{FB}) are rational.

Both the classification of irreducible admissible modules and the
rationality for a vertex operator algebra $V$ are inseparable with
an associative algebra $A(V)$ attached to $V.$ Defined in
\cite{Z}, the associative algebra $A(V)$ which is a quotient of
$V$ plays a fundamental role in the classification of irreducible
admissible modules. It is proved in \cite{Z} and \cite{DLM2} that
there is a one to one correspondence between the equivalence
classes of irreducible admissible $V$-modules and the equivalence
classes of simple $A(V)$-modules. The correspondence is given by
sending the admissible module to its bottom homogeneous subspace
as the admissible module is truncated from below. So in some
sense, the $A(V)$ controls the bottom level of an admissible
module. The relation between admissible modules for a vertex
operator algebra and their bottom levels for $A(V)$ can be
regarded as an analogue of classical highest weight module theory
in the theory of vertex operator algebras. The algebra $A(V)$ is
computable for many vertex operator algebras. The classifications
of irreducible admissible modules for many well known vertex
operator algebras have been achieved using the $A(V)$-theory.
Examples include the affine vertex operator algebras \cite{FZ},
the Virasoro vertex operator algebras \cite{DMZ}, \cite{W1},
lattice type vertex operator algebras \cite{DN1}-\cite{DN3},
\cite{AD}, \cite{TY}, certain $W$-algebras \cite{W2}, \cite{DLTYY}
and some other vertex operator algebras \cite{A2}, \cite{Ad},
\cite{KW}, \cite{KMY}. So the classification of  irreducible
admissible modules for a vertex operator algebra, in principle,
can be done.

It is proved in  \cite{Z} and \cite{DLM2} that if $V$ is rational
then $A(V)$ is a finite dimensional semisimple associative
algebra. A natural question is
$${\rm Does\ the\ semisimplicity\ of\ } A(V)\ {\rm imply\ the\ rationality\
of\ } V?$$ In this paper we give a positive answer to the question
and prove that a simple vertex operator algebra $V$ is rational if
and only if $A(V)$ is semisimple and each irreducible admissible
$V$-module is ordinary.  According to the definition of
rationality, one needs to verify the complete reducibility of any
admissible module to prove the rationality. In practice, this is a
very difficult task. On the other hand, we need the algebra $A(V)$
to classify the irreducible admissible modules. So the $A(V)$ is
available already. While  rationality is an external
characterization of $V$, the semisimplicity of $A(V)$ is certainly
an internal condition on $V$ as the semisimplicity of an
associative algebra can be determined by studying the Jacobson
radical.  We should point out that the assumption that each
irreducible admissible module is ordinary is not  strong as this
is true for all known simple vertex operator algebras.

The main ideas and tools behind the proof of the main result are
the associative algebras $A_n(V)$ \cite{DLM3} and their bimodules
$A_{n,m}(V)$ \cite{DJ1} for nonnegative integers $m,n.$ The
associative algebras $A_n(V)$ are generalizations of $A(V)$ such
that $A_0(V)=A(V).$ We have already  mentioned that $A(V)$
controls the bottom level of an admissible module. The associative
algebra $A_n(V)$,  for each $n\in\Z_{+}$, takes care of the first
$n+1$-levels. Let $M=\bigoplus_{s=0}^{\infty}M(s)$ be an
admissible $V$-module with $M(0)\ne 0$ (see Section 2 for a
precise definition). Then each $M(s)$ for $s\leq n$ is a module
for $A_n(V)$ \cite{DLM3}. Moreover $V$ is rational if and only if
$A_n(V)$ are semisimple for all nonnegative integers $n.$
Theoretically this is a very good result on rationality. But it is
very hard to compute $A_n(V)$ for $n>0$ in practice. Nevertheless,
the $A_n(V)$ theory gives a bridge between a vertex operator
algebra and associative algebras as far as rationality concerns.

Motivated by the fact that $\Hom_{\C}(M(m),M(n))$ is an
$A_n(V)$-$A_m(V)$-bimodule, an abstract $A_n(V)$-$A_m(V)$-bimodule
$A_{n,m}(V)$ is introduced in \cite{DJ1} so that there is a
canonical bimodule homomorphism from  $A_{n,m}(V)$ to
$\Hom_{\C}(M(m),M(n))$ for any admissible $V$-module $M.$ The most
important result about the bimodule theory is an explicit
construction of the Verma type admissible $V$-module generated by
any $A_m(V)$-module $U$ given by
$$M(U)=\bigoplus_{n\geq 0}A_{n,m}(V)\otimes_{A_m(V)}U$$
(see \cite{DJ1} and \cite{DLM3}). As in the classical highest
weight module theory of Lie algebras, $M(U)$ has a unique
irreducible quotient $L(U)$ in case $U$ is an irreducible
$A(V)$-module. The main idea is to prove that $M(U)=L(U)$ if
$A(V)$ is a finite dimensional semisimple associative algebra.

Based on the construction of the Verma type admissible module, a
contravariant pairing between $M(U)$ and $M(U^*)$ is constructed
for any $A_{m}(V)$-module $U$ in this paper.  This pairing is an
analogue of classical contravariant forms for Kac-Moody Lie
algebras and the Virasoro algebra and is used in the proof of the
main theorem. In particular, if $U$ is irreducible, then the left
radical of the pairing is precisely the maximal proper submodule
of $M(U).$ So the contravarint pairing should play an important
role in the study of representation theory for vertex operator
algebras. It is worthy to mention that this contravariant pairing
is totally different from the invariant pairing defined in
\cite{FHL} between any admissible module and its graded dual.

As an application of our main result we prove in this paper that
the orbifold vertex operator algebra $V_L^+$ (see \cite{FLM},
\cite{AD}) is rational for any positive definite even lattice $L.$
The irreducible admissible modules for $V_L^+$ have been
classified in \cite{DN2} and \cite{AD} and each irreducible
admissible module is ordinary. We prove the rationality of $V_L^+$
by showing that $A(V_L^+)$ is semisimple. In \cite{AD}, $A(V_L^+)$
has been understood well enough to classify the irreducible
$A(V_L^+)$-modules. We use a lot of results from \cite{AD} in this
paper to prove the semisimplicity of $A(V_L^+)$ and we refer the
reader to \cite{AD} for a lot of details. If the rank of $L$
equals  one, the rationality of $V_L^+$ has been obtained
\cite{A1} by a different method. Even in this case, it was a very
difficult theorem in \cite{A1}.

It is hard to avoid the regularity \cite{DLM1} and
$C_2$-cofiniteness \cite{Z} when dealing with rationality. It has
been conjectured that the rationality and regularity are
equivalent. There is some progress in proving this conjecture. It
is shown in \cite{L} and \cite{ABD} that $V$ is regular if and
only if $V$ is rational and $C_2$-cofinite. Since any finitely
generated admissible module for a $C_2$-cofinite vertex operator
algebra is ordinary (see \cite{KL} and \cite{ABD}), an immediate
corollary of our result is that $V$ is regular if and only if
$A(V)$ is semisimple and $C_2$-cofinite. Although only the
semisimplicity of $A(V)$ (not the rationality of $V$)  and
$C_2$-cofiniteness were used to obtain the modular invariance of
trace functions \cite{Z} and \cite{DLM5}, the rationality of $V,$
in fact, has already been used by our main result.

It is our hope to remove the assumption that each irreducible
admissible $V$-module is ordinary in the near future.

The main results in this paper had been announced in \cite{DJ2}.

\section{$A_{n,m}(V)$-theory}
\def\theequation{2.\arabic{equation}}
\setcounter{equation}{0}

Let $V=(V,Y,{\bf 1},\omega)$ be a vertex operator algebra (see
\cite{B}, \cite{FLM}).  A {\em weak $V$  module}  is a pair
$(M,Y_M)$, where $M$ is a vector space and $Y_M$ is a linear map
from $V$ to $(\End M)[[z,z^{-1}]]$ satisfying the following axioms
for $u,v\in V$, $w\in M$:
\begin{eqnarray*}\label{e2.1}
& &v_nw=0\ \ \
\mbox{for}\ \ \ n\in \Z \ \ \mbox{sufficiently\ large};\label{vlw0}\\
& &Y_M({\bf 1},z)={\rm id}_{M};\label{vacuum}
\end{eqnarray*}
$$\begin{array}{c}
\displaystyle{z^{-1}_0\delta\left(\frac{z_1-z_2}{z_0}\right)
Y_M(u,z_1)Y_M(v,z_2)-z^{-1}_0\delta\left(\frac{z_2-z_1}{-z_0}\right)
Y_M(v,z_2)Y_M(u,z_1)}\\
\displaystyle{=z_2^{-1}\delta\left(\frac{z_1-z_0}{z_2}\right)
Y_M(Y(u,z_0)v,z_2)}.
\end{array}$$
This completes the definition. We denote this module by $(M,Y_M)$.
An {\em ordinary} $V$-module is a $\C$-graded  weak $V$-module
$$M=\bigoplus_{\lambda \in{\C}}M_{\lambda} $$
such that $\dim M_{\lambda}$ is finite and $M_{\lambda+n}=0$ for
fixed $\lambda$ and $n\in {\Z}$ small enough, where $M_{\lambda}$
is the $\lambda$-eigenspace for $L(0)$ with eigenvalue $\lambda$
and $Y_M(\omega,z)=\sum_{n\in \Z}L(n)z^{-n-2}.$

 An {\em admissible} $V$-module is
a  weak $V$-module $M$ which carries a
${\Z}_{+}$-grading
$$M=\bigoplus_{n\in {\Z}_{+}}M(n)$$
($\Z_+$ is the set of all nonnegative integers) such that if $r,
m\in {\Z} ,n\in {\Z}_{+}$ and $a\in V_{r}$ then
$a_{m}M(n)\subseteq M(r+n-m-1).$ Since the uniform degree shift
gives an isomorphic admissible module we assume $M(0)\ne 0$ in
many occasions.  It is easy to prove that any ordinary module is
an admissible module.

We call a vertex operator algebra {\em rational} if any admissible
module is a direct sum of irreducible admissible modules. It is proved
in \cite{DLM2} that if $V$ is rational then there are only finitely many
irreducible admissible modules up to isomorphism and each irreducible
admissible module is ordinary.

A vertex operator algebra is called {\em regular} if any weak
module is a direct sum of irreducible ordinary modules. It is evident that a
regular vertex operator algebra is rational.

We now review the $A_n(V)$ theory from \cite{DLM2}-\cite{DLM3}. Let $V$ be
a vertex operator algebra. We define two linear operations $*_n$ and $\circ_n$
on $V$ in the following way:
$$u*_nv=\sum_{m=0}^{n}(-1)^m{m+n\choose n}\Res_zY(u,z)\frac{(1+z)^{\wt\,u+n}}{z^{n+m+1}}v$$
$$u\circ_n v=\Res_{z}Y(u,z)v\frac{(1+z)^{\wt u+n}}{z^{2n+2}}
$$
for homogeneous $u,v\in V.$ Let $O_n(V)$ be the linear span of all
$u\circ_n v$ and $L(-1)u+L(0)u.$ Define the linear space $A_n(V)$
to be the quotient $V/O_{n}(V).$

Let $M$ be an admissible $V$-module. For each homogeneous $v\in
V,$ we set $o(v)=v_{\wt v-1}$  on $M$ and extend linearly to whole
$V.$ The main results on $A_n(V)$ were obtained in \cite{DLM3} (
see also \cite{Z} and \cite{DLM2}).
\begin{theorem}\label{tha} Let $V$ be a vertex operator algebra and $n$ a nonnegative
integer. Then

(1) $A_n(V)$ is an associative algebra whose product is induced by
$*_n.$

(2) The identity map on $V$ induces an algebra epimorphism from
$A_n(V)$ to $A_{n-1}(V).$

(3) Let $W$ be a weak module and set
$$\Omega_n(W)=\{w\in W|u_mw=0, u\in V, m> \wt u-1+n\}.$$
Then $\Omega_n(W)$ is an $A_n(V)$-module such that $v+O_n(V)$ acts as $o(v).$

(4) Let $M=\bigoplus_{m=0}^{\infty}M(m)$ be an admissible $V$-module.
Then each $M(m)$ for $m\leq n$ is
an $A_n(V)$-submodule of $\Omega_n(M).$ Furthermore,
$M$ is irreducible if and only if each $M(n)$ is an irreducible
$A_n(V)$-module.

(5) For any $A_n(V)$-module $U$ which is not an
$A_{n-1}(V)$-module, there is a unique Verma type admissible
$V$-module $M(U)$ generated by $U$ so that $M(U)(0)\ne 0$ and
$M(U)(n)=U.$ Moreover, for any weak $V$-module $W$ and any
$A_n(V)$-module homomorphism $f$ from $U$ to $\Omega_n(W)$ there
is a unique $V$-homomorphism from $M(U)$ to $W$ which extends $f.$

(6) $V$ is rational if and only if $A_n(V)$ are finite dimensional
semisimple associative algebras for all $n\geq 0.$

(7) If $V$ is rational then there are only finitely many
irreducible admissible $V$-modules up to isomorphism and each
irreducible admissible module is ordinary.

(8) The linear map  $\phi$ from $V$ to $V$ defined by
$\phi(u)=e^{L(1)}(-1)^{L(0)}u$  for $u\in V$ induces
an anti-involution of $A_n(V).$
\end{theorem}

The algebra $A(V)=A_0(V)$ has been introduced and studied
extensively in \cite{Z} where $*_0$ and $\circ_0$ were denoted by
$*$ and $\circ.$
 Theorem \ref{tha} (1), (3), (5) for $n=0$ have previously been achieved in
\cite{Z} and are very useful in the classification of irreducible
modules for vertex operator algebras  (see \cite{FZ}, \cite{W1},
\cite{W2}, \cite{KW}, \cite{DN1}-\cite{DN3}, \cite{Ad}, \cite{AD},
\cite{DLTYY}, \cite{TY}, \cite{KMY}, \cite{A2}).  We should also
remark that the results in (7) which were parts of the conditions
in the definition of rationality given in \cite{Z} were obtained
in \cite{DLM2}. This shows that the rationality defined in this
paper and \cite{DLM2} is the same as that defined in \cite{Z}.

Next we move to the bimodule theory developed in \cite{DJ1}. For
homogeneous $u\in V,$ $v\in V$ and $m,n,p\in{\mathbb Z}_{+}$,
define the product $\ast_{m,p}^{n}$  on $V$ as follows
$$
u\ast_{m,p}^{n}v=\sum\limits_{i=0}^{p}(-1)^{i}{m+n-p+i\choose
i}{\rm Res}_{z}\frac{(1+z)^{\wt u+m}}{z^{m+n-p+i+1}}Y(u,z)v.
$$
To explain the representation-theoretical meaning of the product
$u\ast_{m,p}^{n}v$ we consider an admissible $V$-module
$M=\bigoplus_{s\geq 0}M(s).$ For homogeneous $u\in V$ we set
$o_t(u)=u_{\wt u-1-t}$ for $t\in \Z.$ Then $o_t(u)M(s)\subset
M(s+t).$ It is proved in \cite{DJ1} that
$o_{n-p}(u)o_{p-m}(v)=o_{n-m}(u\ast_{m,p}^{n}v)$ acting on $M(m).$
If $n=p$, we denote $\ast_{m,p}^{n}$ by $\bar{\ast}_{m}^{n}$, and
if $m=p$, we denote $\ast_{m,p}^{n}$ by $\ast_{m}^{n}$, i.e.,
$$
u\ast_{m}^{n}v=\sum\limits_{i=0}^{m}(-1)^{i}{ n+i\choose i}{\rm
Res}_{z}\frac{(1+z)^{\wt u+m}}{z^{n+i+1}}Y(u,z)v,
$$$$
u\bar{\ast}_{m}^{n}v=\sum\limits_{i=0}^{n}(-1)^{i}{ m+i\choose
i}{\rm Res}_{z}\frac{(1+z)^{\wt u+m}}{z^{m+i+1}}Y(u,z)v.$$

Define  $O'_{n,m}(V)$ to be the linear span of all $u\circ_{m}^{n}v$ and
 $L(-1)u+(L(0)+m-n)u$, where for homogeneous $u\in V$ and $v\in V$,
$$u\circ_{m}^{n}v={\rm
Res}_{z}\frac{(1+z)^{\wt u+m}}{z^{n+m+2}}Y(u,z)v.$$  Then
$O_{n}(V)=O'_{n,n}(V)$ (see \cite{DLM3} and the discussion above).
Let $O''_{n,m}(V)$ be the linear span of
$u\ast_{m,p_{3}}^{n}((a\ast_{p_{1},p_{2}}^{p_{3}}b){\ast}_{m,p_{1}}^{p_{3}}c-
a{\ast}_{m,p_{2}}^{p_{3}}(b{\ast}_{m,p_{1}}^{p_{2}}c)),$
 for $a,b,c,u\in V, p_{1},p_{2},p_{3}\in{\mathbb
Z}_{+}$, and $O'''_{n,m}(V)=\sum_{p\in{\mathbb
Z}_{+}}(V\ast_{p}^{n}O_{p}(V))\ast_{m,p}^{n}V.$ Set
$$
O_{n,m}(V)=O'_{n,m}(V)+O''_{n,m}(V)+O'''_{n,m}(V),
$$
and $$A_{n,m}(V)=V/O_{n,m}(V).$$

Here are the main results on $A_{n,m}(V)$ obtained in \cite{DJ1}.

\begin{theorem}\label{t2.8} Let $V$ be a vertex operator algebra and $m,n$ nonnegative integers. Then

(1) $A_{n,m}(V)$ is an $A_n(V)$-$A_m(V)$-bimodule such that the
left and right actions of $A_n(V)$ and $A_m(V)$ on $A_{n,m}(V)$
are given by $\bar*_m^n$ and $*^n_m$ respectively.

(2) The linear map $\phi: V\to V$ defined by
$\phi(u)=e^{L(1)}(-1)^{L(0)}u$  for $u\in V$ induces
a linear isomorphism from $A_{n,m}(V)$ to $A_{m,n}(V)$
satisfying the following: $\phi(u\bar*_m^n v)=\phi(v)*^m_n\phi(u)$
and  $\phi(v*_m^n w)=\phi(w)\bar*^m_n\phi(v)$ for
$u\in A_n(V),$ $w\in A_m(V)$ and $v\in A_{n,m}(V).$

(3) Let $l$ be nonnegative integers such that $m-l,n-l\geq 0.$ Then $A_{n-l,m-l}(V)$ is an
$A_n(V)$-$A_m(V)$-bimodule and the identity map on $V$ induces an
epimorphism of $A_n(V)$-$A_m(V)$-bimodules from $A_{n,m}(V)$ to
$A_{n-l,m-l}(V).$

(4)  Define a linear map $\psi$: $A_{n,p}(V)\otimes_{A_p(V)}A_{p,m}(V)\rightarrow
A_{n,m}(V)$ by
$$
\psi(u\otimes v)=u\ast_{m,p}^{n}v,$$ for $u\otimes v\in
A_{n,p}(V)\otimes_{A_p(V)}A_{p,m}(V)$. Then $\psi$ is an
$A_{n}(V)-A_{m}(V)$- bimodule homomorphism from
$A_{n,p}(V)\otimes_{A_p(V)}A_{p,m}(V)$ to $A_{n,m}(V)$.

(5) Let $M=\bigoplus_{s=0}^{\infty}M(s)$ be an admissible
$V$-module.  Then $v+O_{n,m}(V)\mapsto o_{n-m}(v)$ gives an
$A_n(V)-A_m(V)$-bimodule homomorphism from $A_{n,m}(V)$ to
$\Hom_{\C}(M(m),M(n)).$

(6) For any $n\geq 0,$ the $A_n(V)$ and $A_{n,n}(V)$ are the same.

(7) Let $U$ be an $A_{m}(V)$-module which can not factor through
$A_{m-1}(V)$. Then
$$
\bigoplus_{n\in{\mathbb Z}_{+}}A_{n,m}(V)\otimes_{A_{m}(V)} U$$ is
a Verma type admissible $V$-module isomorphic to $M(U)$ given in
Theorem \ref{tha} such that  $M(U)(n)=A_{n,m}(V)\otimes_{A_m(V)}U$
with the following universal property: Let $W$ be a weak
$V$-module, then any $A_{m}(V)$-morphism $f: U\rightarrow
\Omega_{m}(W)$ can be extended uniquely to a $V$-homomorphism $f:
M(U)\rightarrow W$.

(8) If $V$ is rational and $W^{j}\!=\!\bigoplus_{n\geq 0}\!W^j(n)$
with $W^j(0)\neq 0$ for $j=1,2,\cdots, s$ are all the inequivalent
irreducible modules of $V,$ then
$$A_{n,m}(V)\cong \bigoplus_{l=0}^{{\rm
min}\{m,n\}}\left(\bigoplus_{i=1}^{s}{\rm Hom}_{\mathbb
C}(W^i(m-l),W^i(n-l))\right).
$$
\end{theorem}

We need the detailed action of $V$ on the Verma type admissible
module $M(U)$ generated by an $A_m(V)$-module $U$ given in Theorem
\ref{t2.8}. For homogeneous $u\in V$ and $p,n\in {\mathbb Z}$, the
component operator $u_{p}$ of $Y_{M(U)}(u,z)=
\sum_{p\in\Z}u_pz^{-p-1}$ which maps $M(U)(n)$ to $M(U)(n+\wt
u-p-1)$ is defined by
\begin{eqnarray}\label{action}
u_{p}(v\otimes w)=\left\{\begin{array}{rl}(u\ast_{m,n}^{\wt
u-p-1+n}v)\otimes w, \quad {\rm if} \ \wt u-1-p+n\geq 0, \\  0,
\quad {\rm if} \ \wt u-1-p+n< 0,
\end{array}\right.
\end{eqnarray}
for $v\in A_{n,m}(V)$ and $w\in U$ (see \cite{DJ1} for details).

Theoretically, the construction of the  Verma type admissible
module $M(U)$ generated by an $A_m(V)$-module $U$ in \cite{DLM3}
is good enough for many purposes.  But $M(U)$ constructed in
[DLM3] is a quotient module for a certain Lie algebra (see
\cite{DLM3}) so it is hard to understand the structure. On the
other hand,   the construction of $M(U)$ given in Theorem
\ref{t2.8} (7) is explicit and each homogeneous subspace $M(U)(n)$
is computable. This new construction of $M(U)$ is fundamental in
our study of rationality in this paper.

\section{Verma type modules}
\def\theequation{3.\arabic{equation}}
\setcounter{equation}{0}

In this section we give a foundational result for this paper. That
is, if $A(V)$ is semisimple for a simple vertex operator algebra
$V$ then the Verma type admissible $V$-module generated by any
irreducible $A(V)$-module is irreducible. We need several lemmas.

Recall that if $A$ is an associative algebra and $U$ a left $A$-module,
then the linear dual $U^*=\Hom_{\C}(U,\C)$ is naturally a right
$A$-module such that $(fa)(u)=f(au)$ for $a\in A,$ $f\in U^*$ and
$u\in U.$
\begin{lem}\label{l2.3} Let $V$ be a vertex operator algebra.
Assume that $A(V)$ is semisimple and $U^{i}$ for $i=1,\cdots,s$ are all the inequivalent
irreducible $A(V)$-modules. Let
${M}(U^{i})=\bigoplus_{n\in{\mathbb Z}_{+}}{M}(U^{i})(n)$ be the
Verma type admissible $V$-module generated by $U^{i}$. Then as an
$A_n(V)$-$A(V)$-bimodule,
$$
A_{n,0}(V)\cong
\bigoplus_{i=1}^{s}{M}(U^{i})(n)\otimes(U^{i})^{*}.$$
\end{lem}

\pf Since $A(V)$ is semisimple, $A_{n,0}(V)$ is a completely
reducible right $A(V)$-module. Note that
$\{(U^1)^*,\cdots,(U^s)^*\}$ is a complete list of inequivalent
irreducible right $A(V)$-modules. Let
$W^i=\Hom_{A(V)}((U^i)^*,A_{n,0}(V))$ which is the multiplicity of
$(U^i)^*$ in $A_{n,0}(V).$ Then $W^i$ is a left $A_n(V)$-module
such that $(af)(x)=af(x)$ for $a\in A_n(V),$ $f\in
\Hom_{A(V)}((U^i)^*,A_{n,0}(V))$ and $x\in (U^i)^*$ and  as
$A_n(V)$-$A(V)$-bimodules
$$ A_{n,0}(V) \cong\bigoplus_{i=1}^sW^{i}\otimes(U^{i})^{*}.$$
To prove the lemma we need to prove that $W^{i}={M}(U^{i})(n)$,
$i=1,\cdots,s$. Recall from  Theorem \ref{t2.8} (7) that for each
$1\leq i\leq s$, ${M}(U^{i})(n)$ is an  $A_{n}(V)$-module and
$${M}(U^{i})=\bigoplus_{n\in{\mathbb
Z}_{+}}A_{n,0}(V)\otimes_{A(V)}U^{i}.$$
Note that $(U^i)^*\otimes_{A(V)}U^j=\delta_{i,j}\C.$ We see immediately
that $W^{i}$ and ${M}(U^{i})(n)$ are isomorphic as required. \qed

For $m,n,p\in{\mathbb Z}_{+}$, let
$$A_{n,p}(V)\ast_{m,p}^{n}A_{p,m}(V)=\{a\ast_{m,p}^{n}b|a\in
A_{n,p}(V),b\in A_{p,m}(V)\}.$$ Then
$A_{n,p}(V)\ast_{m,p}^{n}A_{p,m}(V)$ is an
$A_{n}(V)-A_{m}(V)$-subbimodule of $A_{n,m}(V)$ by Theorem \ref{t2.8}
(4). Also recall from Theorem \ref{t2.8} (4) the  $A_{n}(V)-A_{m}(V)$-bimodule homomorphism $\psi:$
 $A_{n,p}(V)\otimes_{A_{p}(V)}A_{p,m}(V)\rightarrow A_{n,m}(V)$ defined  by
 $$
\psi(u\otimes v)=u\ast_{m,p}^{n}v,$$ for $u\in A_{n,p}(V)$ and
$v\in A_{p,m}(V)$.

\begin{lem}\label{l2.11} Let $V$ be a simple vertex operator algebra
and  $A(V)$ be  semisimple. Then for  integer $m\geq 2$, we
have  $$A_{0,m}(V)\ast_{0,m}^{0}A_{m,0}(V)=A(V).$$
\end{lem}

\pf It is easy to see that $A_{0,n}(V)\ast_{0,n}^{0}A_{n,0}(V)$ is
a two-sided ideal of $A(V)$. Let $U$ be an irreducible module of
$A(V)$ and suppose that for some positive integer $n$,
$$(A_{0,n}(V)\ast_{0,n}^{0}A_{n,0}(V))\otimes_{A(V)} U=0.$$
 Let $M(U)$ be the Verma type admissible $V$-module generated by $U$ and $X$ the
 admissible $V$-submodule of $M(U)$ generated by $M(U)(n)$.
Then by  Proposition 4.5.6 of \cite{LL} (see also \cite{DM1}), $X$
is spanned by $u_pM(u)(n)$ for $u\in V$ and $p\in\Z.$ In
particular, $X(0)$ is spanned by $u_{\wt u-1+n}M(u)(n)$ for $u\in
V.$ From the module construction of $M(U)$ given in Theorem
\ref{t2.8} and the action of $V$ on $M(U)$ (\ref{action}) we see
that
$$X(0)=(A_{0,n}(V)\ast_{0,n}^{0}A_{n,0}(V))\otimes_{A(V)} U=0.$$
So $X$ is a proper admissible submodule of $M(U).$

Clearly, for any non-zero element $u$ in $U$, we have $L(1)u=0$
and  $L(-n)u\in M(U)(n)$ from the construction of $M(U).$ Thus
$0=L(1)^nL(-n)u\in X(0).$ It follows that $L(0)u=0$. This shows
that  $L(0)U=0$.

Let $W(U)=M(U)/X$ be the quotient module, then $W(U)(n)=0$ as
$M(U)(n)\subset X.$ Since $W(U)(0)=U$, we can assume that
$W(U)(n-k)\neq 0$, $W(U)(n-k+1)=0$, for some positive integer $k$.
Let $v$ be a non-zero element in $W(U)(n-k)$, then $L(-1)v=0$. By
Corollary 4.7.6 and Proposition 4.7.9 of \cite{LL}, $v$ is a
vacuum-like vector  of $W(U)$ and the admissible $V$-submodule of
$W(U)$ generated by $v$ is isomorphic to $V$. In particular,
$L(0)v=0$. On the other hand, $L(0)v=(n-k)v$ as $L(0)U=0$. This
deduces that $k=n$. Since $V_2$ contains the Virasoro element
$\omega\ne 0$ and $L(-1)$ from $V_m$ to $V_{m+1}$ is injective for
any $m\geq 1,$ we assert that $V_{m}\neq 0$ if $m\geq 2$. This
implies that $n=1$. So for any $m\geq 2$, we have
$$(A_{0,m}(V)\ast_{0,m}^{0}A_{m,0}(V))\otimes U\neq 0.$$
Since $A(V)$ is a finite dimensional semisimple associative
algebra, the lemma follows. \qed

\begin{lem} \label{l2.6} Let $V$
be a  simple vertex operator algebra such that  $A(V)$ is
semisimple. Let $n\in{\mathbb Z}_{+}$ greater than 1, then the
$A(V)-A(V)$-bimodule homomorphism $\psi$ from
$A_{0,n}(V)\otimes_{A_{n}(V)}A_{n,0}(V)$ to $A(V)$ is an
isomorphism.
\end{lem}

\pf  By Lemma \ref{l2.11}, $A_{0,n}(V)\ast_{0,n}^{0}A_{n,0}(V)=A(V)$. So it is enough to prove that
 $\psi$ is injective. Let ${\bf 1}$  be the
identity  of $A(V)$ and $u^{j}\in A_{0,n}(V)$ and $v^{j}\in
A_{n,0}(V)$, $j=1,2,\cdots,k$, be such that
$\sum\limits_{j=1}^{k}u^{j}\ast_{0,n}^{0}v^{j}=0$. That is,
$\sum\limits_{j=1}^{k}u^{j}\otimes_{A_n(V)}v^{j}$ lies in the
kernel of $\psi.$ By Lemma \ref{l2.11},  there exist $a^{i}\in
A_{0,n}(V)$, $b^{i}\in A_{n,0}(V)$, $i=1,2,\cdots,r$ such that
$$\sum\limits_{i=1}^{r}a^{i}\ast_{0,n}^{0}b^{i}={\bf 1}.$$
Then
\begin{eqnarray*}
& & \ \ \ \sum\limits_{j=1}^{k}u^{j}\otimes
v^{j}=(\sum\limits_{j=1}^{k}u^{j}\otimes v^{j})\cdot{\bf
1}=\sum\limits_{j=1}^{k}u^{j}\otimes (v^{j}\ast_{0,0}^{n}{\bf 1})\\
& &=\sum\limits_{j=1}^{k}u^{j}\otimes
(v^{j}\ast_{0,0}^{n}(\sum\limits_{i=1}^{r}a^{i}\ast_{0,n}^{0}b^{i}))
=\sum\limits_{j=1}^{k}\sum\limits_{i=1}^{r}u^{j}\otimes((v^{j}\ast_{n,0}^{n}a^{i})\ast_{0,n}^{n}b^{i})\\&
&
=\sum\limits_{j=1}^{k}\sum\limits_{i=1}^{r}(u^{j}\ast_{n,n}^{0}(v^{j}\ast_{n,0}^{n}a^{i}))\otimes
b^{i}=\sum\limits_{j=1}^{k}\sum\limits_{i=1}^{r}((u^{j}\ast_{0,n}^{0}v^{j})\ast_{n,0}^{0}a^{i})\otimes
b^{i}\\ & &= 0.
\end{eqnarray*}
This means that $\psi$ is injective. \qed

We are now in a position to prove the following important result
on the Verma type admissible module.

\begin{theorem}\label{t2.7} Let $V$
be a simple  vertex operator algebra such that $A(V)$ is
semisimple. Let $U$ be an irreducible module of $A(V)$, then the
Verma type admissible $V$-module $M(U)=\bigoplus_{n\in{\mathbb
Z}_{+}}A_{n,0}(V)\otimes_{A(V)}U$ generated by $U$ is
irreducible.
\end{theorem}

\pf For  integer $n\geq 1$, set
$$
S(n)=\{x\in A_{n,0}(V)|u\ast_{0,n}^{0} x=0, u\in A_{0,n}(V)\}.$$
Then by the fact that
$u\ast_{0,n}^{0}(x\ast_{0,0}^{n}v)=(u\ast_{0,n}^{0}x)\ast v$ and
$u\ast_{0,n}^{0}(a\ast_{0,n}^{n}x)=(u\ast_{n,n}^{0}a)\ast_{0,n}^{0}x$
for all $u\in A_{0,n}(V),$ $x\in A_{n,0}(V),$ $a\in A_{n}(V),$ $v\in
A(V)$, $S(n)$ is an $A_{n}(V)-A(V)$-subbimodule of $A_{n,0}(V)$.
The $S(n)$ can be understood naturally from the Verma type admissible
$V$-module
$$M(A(V))=\bigoplus_{n\geq 0}M(A(V))(n)=\bigoplus_{n\geq 0}A_{n,0}(V)\otimes_{A(V)}A(V)=\bigoplus_{n\geq 0}A_{n,0}(V)$$
generated by $A(V).$ Let $S$ be the maximal submodule of $M(A(V))$
such that
$$S\cap M(A(V))(0)=S\cap A(V)=0.$$
 Then it is clear that
$S=\sum_{n\geq 1}S(n).$

As before, let $U^{i}$ for $i=1,\cdots,s$ be the inequivalent
irreducible modules (which are necessarily finite dimensional) of $A(V)$.
Then $A(V)=\bigoplus_{i=1}^{s}{\rm End}_{\mathbb C}(U^i)$
and $U^i$ can be considered as a
simple left ideal of $A(V)$ and the action of $A(V)$ on $U^i$ is
just the multiplication of $A(V)$. Let
$$M(U^i)=\bigoplus_{n\in{\mathbb
Z}_{+}}M(U^i)(n)=\bigoplus_{n\in{\mathbb
Z}_{+}}A_{n,0}(V)\otimes_{A(V)} U^{i}$$
be the Verma type
 admissible $V$-module generated by $U^i$. Then
$$M(U^i)(n)=A_{n,0}(V)\otimes_{A(V)} U^i=A_{n,0}(V)\ast_{0,0}^{n}U^i\subseteq A_{n,0}(V)$$
where we identify $U^i$ with a simple left ideal of $A(V).$ In
fact, we can regard $M(U^i)$ as an admissible submodule of
$M(A(V)).$ The containment $M(U^i)(n)\subset  A_{n,0}(V)$ can be
understood easily from Lemma \ref{l2.3}.

Let $J(U^i)=\bigoplus_{n\in{\mathbb Z}_{+}}J(U^i)(n)$ be the
maximal proper admissible $V$-submodule of $M(U^i)$. Then the
quotient module $W^i=M(U^i)/J(U^i)=\bigoplus_{n\in{\mathbb
Z}_{+}}M(U^i)(n)/J(U^i)(n)$ is an irreducible admissible
$V$-module. By Theorem \ref{tha} (4), $W^i(n)$ is an irreducible $A_{n}(V)$-module.
Since $U^i$ is an irreducible $A(V)$-module, we have
$J(U^i)(0)=0$. Regarding $M(U^i)$ as an admissible submodule of
$M(A(V))$ then $J(U^i)(n)$ is a subspace of $M(A(V))(n)=
A_{n,0}(V)\otimes_{A(V)}A(V)=A_{n,0}(V)$ and
\begin{equation}\label{ea2.1}
A_{0,n}(V)\ast_{0,n}^{0}J(U^i)(n)=0,\end{equation}
 for all positive
integer $n$. That is, $J(U^i)(n)$ is a subspace of $S(n)$ for all
$i.$ To prove that any Verma type admissible $V$-module
$M(U)=\bigoplus_{n\in{\mathbb Z}_{+}}A_{n,0}(V)\otimes_{A(V)}U$
generated by an irreducible $A(V)$-module $U$ is irreducible, it
is enough to prove that $S(n)=0$ for all $n\geq 1.$ From the
assumption, $V$ has only finitely many irreducible admissible
modules, so $L(-1)$ from $S(n)$ to $S(n+1)$ is injective if $n$ is
large. So it suffices to show that $S(n)=0$ for all large $n.$

{\bf Claim:} If $n\geq 2,$ then $(S(n)\ast_{n,0}^{n}A_{0,n}(V))\otimes_{A_{n}(V)}W^i(n)=0$.

If $(A_{n,0}(V)\ast_{n,0}^{n}A_{0,n}(V))\otimes_{A_{n}(V)}
W^i(n)=0$, then
$(S(n)\ast_{n,0}^{n}A_{0,n}(V))\otimes_{A_{n}(V)}W^i(n)=0$.
So we now  assume
that $(A_{n,0}(V)\ast_{n,0}^{n}A_{0,n}(V))\otimes_{A_{n}(V)} W^i(n)\neq
0.$ By Lemma \ref{l2.6} and (\ref{ea2.1}) we
have for $n\geq 2$ that
$$A_{0,n}(V)\otimes_{A_{n}(V)}J(U^i)(n)=0.$$
This implies that
\begin{equation}\label{ea2.2}
A_{0,n}(V)\otimes_{A_{n}(V)}M(U^i)(n)\cong
A_{0,n}(V)\otimes_{A_{n}(V)}W^i(n).\end{equation} Consider the
Verma type admissible $V$-module
$M(W^i(n))=\bigoplus_{m\in{\mathbb
Z}_{+}}A_{m,n}(V)\otimes_{A_{n}(V)} W^i(n)$ generated by $W^i(n)$.
Then $M(W^i(n))(n)=W^i(n)$ is an irreducible $A_n(V)$-module and
$(A_{n,0}(V)\ast_{n,0}^{n}A_{0,n}(V))\otimes_{A_{n}(V)} W^i(n)$ is
a nonzero submodule of $W^i(n).$ This forces
$(A_{n,0}(V)\ast_{n,0}^{n}A_{0,n}(V))\otimes_{A_{n}(V)}
W^i(n)=W^i(n).$
From  the definition of $S(n)$, we have
$$[A_{0,n}(V)\ast_{n,n}^{0}(S(n)\ast_{n,0}^{n}A_{0,n}(V))]\otimes_{A_{n}(V)}
 W^i(n)$$
 $$=[(A_{0,n}(V)\ast_{0,n}^{0}S(n))\ast_{n,0}^{0}A_{0,n}(V)]\otimes_{A_{n}(V)}W^i(n)=0.$$
 So $(S(n)\ast_{n,0}^{n}A_{0,n}(V))\otimes_{A_{n}(V)}
 W^i(n)$ is a proper $A_{n}(V)$-submodule of
irreducible $A_n(V)$-module
$(A_{n,0}(V)\ast_{n,0}^{n}A_{0,n}(V))\otimes_{A_{n}(V)}
 W^i(n).$ Thus
$$(S(n)\ast_{n,0}^{n}A_{0,n}(V))\otimes_{A_{n}(V)}
 W^i(n)=0.$$
This establishes the claim.

By (\ref{ea2.2}), we have
$$(S(n)\ast_{n,0}^{n}A_{0,n}(V))\otimes_{A_{n}(V)}
 M(U^i)(n)=0, \ i=1,2,\cdots,s.$$
Applying Lemma \ref{l2.3} yields
\begin{eqnarray*}
& & \ \ \
(S(n)\ast_{n,0}^{n}A_{0,n}(V))\otimes_{A_{n}(V)}A_{n,0}(V)\\& &=
(S(n)\ast_{n,0}^{n}A_{0,n}(V))\otimes_{A_{n}(V)}(\bigoplus_{i=1}^{s}
 M(U^i)(n)\otimes (U^{i})^*)\\
& &=0.
 \end{eqnarray*}
On the other hand, by Lemma \ref{l2.11} we have
$$
S(n)=S(n)\ast_{0,0}^{n}A(V)=S(n)\ast_{0,0}^{n}(A_{0,n}(V)\ast_{0,n}^{0}A_{n,0}(V))$$
which is exactly
$(S(n)\ast_{n,0}^{n}A_{0,n}(V))\otimes_{A_{n}(V)}A_{n,0}(V)$ by
Lemma \ref{l2.6}. As a result, $S(n)=0.$ This finishes the proof.
\qed

We should point out that if $A(V)$ is not semisimple, Theorem
\ref{t2.7} is false. Here is a counter example. Recall that the
abstract Virasoro algebra $\Vir$ has a basis $\{L_n|n\in\Z\}\cup
\{c\}$ with relation:
$$[L_m,L_n]=(m-n)L_{m+n}+\frac{m^3-m}{12}\delta_{m+n,0}c$$
and $c$ is in the center. Then
$${\rm Vir}^{\geq -1}=\oplus_{n=-1}^{\infty}\C L_n\oplus \C c$$
is a subalgebra. Given a complex number $k,$ consider the induced
module
$$V(k)=U(\Vir)\otimes _{U(\Vir^{\geq -1})}\C_{k}$$
where $\C_{k}=\C$ is a module for $\Vir^{\geq -1}$ such that
$L_n1=0$ for $n\geq -1$ and $c1=k.$ Then $V(k)$ is a vertex
operator algebra (see \cite{FZ}). If $k=1,$ then $V(1)$ is a
simple vertex operator algebra as $V(1)$ is an irreducible module
for the Virasoro algebra. It is computed in \cite{FZ} that
$A(V(1))$ is isomorphic to the polynomial algebra $\C[x]$ which is
not semisimple. Let $U$ be an irreducible $A(V(1))$-module such
that $\omega$ acts as a positive integer $m.$ Then the Verma type
admissible $V(1)$-module $M(U)$ generated by $U$ is the Verma
module $V(1,m)= U(Vir)\otimes_{U(\Vir^{\geq 0})}\C_{1,m}$ for
\Vir\ where $\Vir^{\geq 0}= \sum_{n=0}^{\infty}\C L_n\oplus\C c$
acts on $\C_{1,m}=\C$ in the following way: $L_n1=0$ for $n>0$ and
$L_01=m,\ c1=1.$ Clearly, $V(1,m)$ is not irreducible (see
\cite{KR}).

One can also find counter example from the affine vertex operator algebra
if $A(V)$ is not semisimple.

\section{Bilinear pairings}
\def\theequation{4.\arabic{equation}}
\setcounter{equation}{0}

Let $M=\bigoplus_{n\geq 0}M(n)$ be an admissible $V$-module. Then
there is a $V$-invariant bilinear pairing $(\cdot,\cdot)$ between
$M$ and its contragradient module $M'=\sum_{n\geq 0}M(n)^*$ in the
sense that
$$(Y(u,z)w',w)=(w',Y(e^{zL(1)}(-z^{-2})^{L(0)}u,z^{-1})w)$$
for $u\in V,$ $w\in M$ and $w'\in M'$ where
$M(n)^*=\Hom_{\C}(M(n),\C)$ (see \cite{FHL}). The construction of
contragradient module $M'$ and the non-degenerate pairing are very
useful in the theory of vertex operator algebras. But this
bilinear pairing is different from the contravariant form on the
Verma modules for affine Lie algebras or the Virasoro algebra when
$V$ is an affine or Viraoro vertex operator algebra. While the
contravariant form on a Verma module in the classical case is
degenerate in general,
 the invariant bilinear pairing defined in \cite{FHL} is always
non-degenerate. In this section we will construct a different
invariant bilinear pairing between the two Verma type admissible
$V$-modules $M(U)$ and $M(U^*)$ for any $A_{m}(V)$-module $U$,
where $m\in\Z_{+}$. This pairing is an analogue of classical
contravariant forms for Kac-Moody Lie algebras and the Virasoro
algebra. This construction heavily depends on the construction of
$M(U)$ given in Theorem \ref{t2.8} in terms of bimodules.

For $m\in\Z_{+}$, let $U$ be an $A_{m}(V)$-module and recall the
anti-involution $\phi$ of $A_{m}(V)$ from Theorem \ref{tha}. Then
the dual space $U^*$ of $U$ is an $A_{m}(V)$-module under the
following action:
\begin{equation}\label{ea2.3}
(u\cdot f)(x)=f(\phi(u)\cdot x)=(f, \phi(u)\cdot x),\end{equation}
for $u\in A_{m}(V)$, $f\in U^*$ and $x\in U$.

Let $M(U)=\bigoplus_{n\in{\mathbb
Z}_{+}}A_{n,m}(V)\otimes_{A_{m}(V)} U$ and
$M(U^*)=\bigoplus_{n\in{\mathbb
Z}_{+}}A_{n,m}(V)\otimes_{A_{m}(V)} U^*$ be the Verma type
admissible $V$-modules generated by $U$ and $U^*$ respectively.
Recall from Theorem \ref{t2.8} the linear isomorphism $\phi$ from
$A_{n,m}(V)$ to $A_{m,n}(V).$ We define a bilinear pairing
$(\cdot,\cdot)$ on $M(U^*)\times M(U)$ as follows:
$$
(x\otimes f, y\otimes u)=(f, [(\phi(x)\ast_{m,n}^{m} y)]\cdot
u),$$ for $x,y\in A_{n,m}(V), f\in U^*, u\in U,n\in{\mathbb
Z}_{+}$; and
$$
(A_{p,m}(V)\otimes_{A_{m}(V)}U^*,A_{n,m}(V)\otimes_{A_{m}(V)}U)=0$$
 for $p\neq n$. That is $(M(U^*)(p),M(U)(n))=0$ if $p\ne n.$

\begin{lem} The bilinear pairing  $(\cdot,\cdot)$ is well defined.
\end{lem}

\pf  Let
 $x\in A_{n,m}(V), y\in A_{p,m}(V), f\in U^*, v\in
U,p,n\in{\mathbb Z}_{+}$ and $a\in A_{m}(V).$  If $p=n$, we have
from Theorem \ref{t2.8}  that
\begin{eqnarray*}
& & ((x\ast_{m,m}^{n}a)\otimes f,y\otimes v)=(f,[\phi(x\ast_{m,m}^{n}a)\ast_{m,n}^{m}y]\cdot v)\\
& &\ \ \ \
=(f,[(\phi(a)\ast_{n,m}^{m}\phi(x))\ast_{m,n}^{m}y]\cdot
v)\\
& &\ \ \ \ =(f,[\phi(a)\ast_{m}(\phi(x)\ast_{m,n}^{m}y)]\cdot v)\\
& &\ \ \ \ =(f,\phi(a)\cdot[(\phi(x)\ast_{m,n}^{m}y)\cdot v])\\
& &\ \ \ \ =(a\cdot
f,(\phi(x)\ast_{m,n}^{m}y)\cdot v)\\
& &\ \ \ \ =(x\otimes (a\cdot f), y\otimes v)
\end{eqnarray*}
where we have used (\ref{ea2.3}) in the fifth equality.

If $p\neq n$, it is clear that
$$((x\ast_{m,m}^{n}a)\otimes f,y\otimes
v)=(x\otimes (a\cdot f), y\otimes v)=0.$$
 Similarly, we have
$$
(x\otimes f, (y\ast_{m,m}^{p}a)\otimes v)=(x\otimes f, y\otimes
(a\cdot v)).$$ The proof is complete. \qed

\begin{prop} \label{p2.3} The bilinear pairing
$(\cdot,\cdot)$ on $M(U^*)\times M(U)$ is
invariant in the sense that
\begin{equation}\label{4.2}
(Y_{M}(u,z)(x\otimes f),y\otimes v)=(x\otimes f,
Y_{M}(e^{zL(1)}(-z^{-2})^{L(0)}u,z^{-1})(y\otimes v))
\end{equation}
for $x\in A_{n,m}(V), y\in A_{p,m}(V), f\in U^*, v\in
U,n,p\in{\mathbb Z}_{+}$ and $u\in V$.
\end{prop}

\pf It is enough to  prove the coefficients of $z^{-q-1}$ in both
sides of (\ref{4.2}) are equal for all  $q\in{\mathbb Z}.$ That
is, we only need to prove
\begin{equation}\label{ea0}
(u_{q}(x\otimes f),y\otimes v)=(x\otimes f,
\sum\limits_{j=0}^{\infty}\frac{(-1)^{\wt
u}}{j!}(L(1)^{j}u)_{-q+2\wt u-j-2}(y\otimes v)),
\end{equation}
for $x\in A_{n,m}(V),y\in A_{p,m}(V), f\in U^*, v\in
U,n,p\in{\mathbb Z}_{+}, q\in{\mathbb Z}$ and homogeneous $u\in
V$.

First assume that ${\wt}u-q-1+n\neq p$. Since $u_{q}
M(U^*)(n)\subset M(U^*)({\wt u}-q-1+n),$ then
$(u_qM(U^*)(n),M(U)(p))=0$ from the definition of $(\cdot,\cdot).$
Similarly,
$$(M(U^*)(n),\sum\limits_{j=0}^{\infty}\frac{(-1)^{{\wt }u}}{j!}(L(1)^{j}u)_{-q+2{\wt }u-j-2}M(U)(p))=0.$$
So (\ref {ea0}) is true if ${\wt }u-q-1+n\neq p$.

If ${\wt}u-q-1+n=p$, using Theorem \ref{t2.8}  and
 the action of $u_q$ given in (\ref{action}) we have
\begin{eqnarray*}
& & \ \ \ (u_{q}(x\otimes f),y\otimes v)
=((u\ast_{m,n}^{p}x)\otimes f,y\otimes v)\\
& & =(f,(\phi(u\ast_{m,n}^{p}x)\ast_{m,p}^{m}y)\cdot
v)=(f,((\phi(x)\ast_{p,n}^{m}\phi(u))\ast_{m,p}^{m}y)\cdot v)\\
& &=(f,(\phi(x)\ast_{m,n}^{m}(\phi(u)\ast_{m,p}^{n}y))\cdot
v)=(x\otimes f,(\phi(u)\ast_{m,p}^{n}y)\otimes v)\\
& & =(x\otimes f,[(\sum\limits_{j=0}^{\infty}\frac{(-1)^{{\wt
}u}}{j!}L(1)^{j}u)\ast_{m,p}^{n}y]\otimes
v)\\
& &=(x\otimes f,(\sum\limits_{j=0}^{\infty}\frac{(-1)^{{\wt
}u}}{j!}L(1)^{j}u)_{-q+2{\wt }u-j-2}(y\otimes v)),
\end{eqnarray*}
as desired. \qed

The proof of the following corollary  is standard.
\begin{coro}\label{t2.4} Let $V$ be a  vertex operator
algebra and $m\in\Z_{+}$. Let  $U$ be an $A_{m}(V)$-module. Then
$$
J(U)=\{w\in M(U)|(w', w)=0, w'\in M(U^*)\}
$$
is the maximal proper admissible $V$-submodule of $M(U)$ such that
$J(U)\cap M(U)(m)=0.$ In particular, if $U$ is irreducible then
$J(U)$ is the unique maximal proper admissible submodule of
$M(U).$
\end{coro}

Corollary \ref{t2.4} tells us that the bilinear pairing defined in
this section is an analogue of the classical contravariant form in
the theory of vertex operator algebras. This result is certainly
an important application of the construction of the Verma type
admissible $V$-module $M(U)=\bigoplus_{n\geq
0}A_{n,m}(V)\otimes_{A_{m}(V)}U$ generated by an $A_{m}(V)$-module
$U.$ It is hard to imagine how to define the invariant bilinear
pairing on $M(U^*)\times M(U)$ without this module construction.

By Theorem \ref{t2.7} and Corollary \ref{t2.4}, we immediately have
\begin{coro}\label{co2.6} Let $V$
be a  simple vertex operator algebra  such that $A(V)$ is
semisimple. Let $U$ be an irreducible $A(V)$-module, then the
bilinear form $(\cdot,\cdot)$ on $M(U^*)\times M(U)$ is
non-degenerate.
\end{coro}

\section{Rationality}
\def\theequation{5.\arabic{equation}}
\setcounter{equation}{0}

We prove in this section that a simple vertex operator algebra $V$
is rational if and only if $A(V)$ is semisimple and each
irreducible admissible $V$-module is ordinary. This is the main
result of this paper.

Recall from Theorem \ref{t2.8} that $A_{m,0}(V)\ast_{m,0}^{m}A_{0,m}(V)
\subset A_m(V).$ We have
\begin{lem}\label{l2.8}
Let $V$ be a vertex operator algebra. Then
$$A_{m,0}(V)\ast_{m,0}^{m}A_{0,m}(V)\cong O_{m-1}(V)/O_{m}(V),$$
for $m\in\Z_{+}\setminus \{0\}$.
\end{lem}

\pf For $m\in\Z_{+}\setminus \{0\}$, let
$M(A_{m}(V))=\bigoplus_{n\in\Z_{+}}A_{n,m}(V)$ be the Verma type
admissible $V$-module generated by $A_{m}(V)$. Then
$$M(A_{m}(V))(0)=A_{0,m}(V).$$
Let $M'$ be the $V$-submodule of $M(A_{m}(V))$ generated by
$A_{0,m}(V)$. Then
$$M'(m)=A_{m,0}(V)\ast_{m,0}^{m}A_{0,m}(V).$$
Set $W=M(A_{m}(V))/M'$. Then $W$ is an admissible $V$-module such
that $W(m)=A_{m}(V)/(A_{m,0}(V)\ast_{m,0}^{m}A_{0,m}(V))$ and
$W(0)=0$. So for all homogeneous $u\in V$, we have
$$u_{{\wt }u-1+m}W(m)=0.$$
This means that $W(m)$ is an $A_{m-1}(V)$-module by Theorem \ref{tha}.
Thus
\begin{equation}\label{ea2.0}
O_{m-1}(V)\ast_{m}A_{m}(V)\subset A_{m,0}(V)\ast_{m,0}^{m}A_{0,m}(V).
\end{equation}
Note that $A_{m}(V)=V/O_m(V)$  and
$A_{m,0}(V)\ast_{m,0}^{m}A_{0,m}(V)$ is the linear span of
$u\ast_{m,0}^{m}v+O_{m}(V)$ for all $v\in V$ and homogeneous $u\in
V$. Let $B_{m}(V)$ be the linear span of $u\ast_{m,0}^{m}v$ for
all $v\in V$ and homogeneous $u\in V$. Then
$$W(m)\cong V/(B_{m}(V)+O_{m}(V)), \ \ A_{m,0}(V)\ast_{m,0}^{m}A_{0,m}(V)=(B_{m}(V)+O_{m}(V))/O_{m}(V).$$
From (\ref{ea2.0}), we have
$$
O_{m-1}(V)\subseteq B_{m}(V)+O_{m}(V).$$ On the other hand, for $v\in V$ and homogeneous $u\in V$
$$
u\ast_{m,0}^{m}v={\rm Res}_{z}\frac{(1+z)^{{\wt
}u+m}}{z^{2m+1}}Y(u,z)v\in O_{m-1}(V).$$ That is, $B_m(V)\subset
O_{m-1}(V).$ As a result, we have
$$B_{m}(V)+O_{m}(V)=O_{m-1}(V),$$
and therefore the lemma holds. \qed

\begin{lem}\label{p2.4} Let $V$
be a  vertex operator algebra such that  $A(V)$ is semisimple. Let
$U^{i}$ for $i=1,2,\cdots,s$ be
 all the inequivalent irreducible
modules of $A(V)$. Then for every positive integer $n$, we have
$$
A_{0,n}(V)\cong\bigoplus_{i=1}^{s}U^{i}\otimes M((U^i)^*)(n),$$
where the left action of $A(V)$ on $(U^i)^*$ is defined by
(\ref{ea2.3}).
\end{lem}

\pf By Lemma \ref{l2.3}, we have
$$
A_{n,0}(V)\cong
\bigoplus_{i=1}^{s}{M}(U^{i})(n)\otimes{(U^{i})}^{*}.$$
$A_{n,0}(V)$ is an $A_{n}(V)-A(V)$-bimodule with the following
left   and right actions by $A_{n}(V)$ and $A(V)$ respectively:
$$
a\cdot(\sum\limits_{i=1}^{s}u^{i}\otimes
v^{i})=\sum\limits_{i=1}^{s}(a\cdot u^{i})\otimes v^{i}, \ \
(\sum\limits_{i=1}^{s}u^{i}\otimes v^{i})\cdot
b=\sum\limits_{i=1}^{s}u^{i}\otimes v^{i}\cdot b,
$$
where $a\in A_{n}(V)$, $b\in A(V)$, $u^{i}\in{M}(U^{i})(n)$,
$v^{i}\in{(U^{i})}^{*}$,  and $(v^{i}\cdot b)(c)=v^{i}(b\cdot c)$,
for $c\in U^i$. By Theorem \ref{t2.7}, each
${M}(U^{i})(n)\otimes{(U^{i})}^{*}$ is an irreducible
$A_{n}(V)-A(V)$-bimodule.

Let $m,n\in\Z_{+}$, recall from Theorem \ref{t2.8} that $\phi$ is
the linear map from $A_{n,m}(V)$ to $A_{m,n}(V)$ defined as
$\phi(u)=e^{L(1)}(-1)^{L(0)}u$, for $u\in A_{n,m}(V)$. By
Proposition 3.2 of \cite{DJ1}, $\phi$, in fact,  is an isomorphism
of $A_{n}(V)-A_{m}(V)$-bimodules from $A_{n,m}(V)$ to
$A_{m,n}(V)$, where the left action $\bar\cdot_m^n$ of $A_n(V)$
and the right action $\cdot_m^n$ of $A_m(V)$ on $A_{m,n}(V)$ are
defined by $u\bar\cdot_m^nv=v*_n^m\phi(u)$ and $v\cdot_m^n
a=\phi(a)\bar *^m_n v$ respectively for $u\in A_n(V),$ $v\in
A_{m,n}(V)$ and $a\in A_m(V).$

So it is sufficient to show that $A_{n,0}(V)$ is isomorphic to
$\bigoplus_{i=1}^{s}U^{i}\otimes M((U^i)^*)(n)$  as
$A(V)-A_n(V)$-bimodule with the new actions. As an
$A(V)-A_n(V)$-bimodule $M(U^i)(n)\otimes (U^i)^*$ is clearly
isomorphic to $(U^i)^*\otimes M(U^i)(n).$ Since each $U^i$ is
finite dimensional we finish the proof. \qed

We now prove the main theorem.
\begin{theorem}\label{t2.11} Let $V$ be a simple vertex operator
algebra. Then $V$ is rational if and only if $A(V)$ is semisimple
and each irreducible admissible $V$-module is ordinary.
\end{theorem}

\pf By Theorem \ref{tha}, if $V$ is rational then $A(V)$ is semisimple
and each irreducible admissible $V$-module is ordinary. Now we assume that
$A(V)$ is semisimple and each irreducible admissible $V$-module is ordinary.
By Theorem \ref{tha},  it is good enough to prove that
$A_{n}(V)$ are semisimple for all positive integers $n.$ We achieve this by
induction on $n.$

Suppose $n\geq 1$ such that $A_{m}(V)$ are all semisimple for
$0\leq m\leq n-1.$ In order to prove that $A_n(V)$ is semisimple,
it is good enough to prove that any indecomposable $A_n(V)$-module
$X$ is irreducible. We may assume that $X$ is not an
$A_{n-1}(V)$-module.

By Lemma \ref{l2.3} and Lemma \ref{p2.4}, both $A_{n,0}(V)$ and
$A_{0,n}(V)$ are finite dimensional (here we are using the
assumption that each irreducible admissible module is ordinary).
It follows from Lemma \ref{l2.8} that $O_{n-1}(V)/O_{n}(V)$ is
finite dimensional. Since $A_{n-1}(V)$ is also finite dimensional
by inductive assumption, we see that $A_n(V)$ is finite
dimensional. In fact, the dimension of $A_n(V)$ is the sum of
dimensions of $O_{n-1}(V)/O_{n}(V)$ and $A_{n-1}(V).$ This implies
that $X$ is finite dimensional.

Consider the Verma type admissible $V$-module
$M(X)=\bigoplus_{k\geq 0}A_{k,n}(V)\otimes_{A_n(V)}X.$ Then by
Theorem \ref{tha}, $M(X)(0)=A_{0,n}(V)\otimes_{A_n(V)}X \ne 0.$ By
assumption, $M(X)(0)$ is a semisimple $A(V)$-module. Let $Z$ be
the submodule of $M(X)$ generated by $M(X)(0).$ Then $Z$ is
completely reducible by Theorem \ref{t2.7} and $Z(n)=X\cap Z$ is a
completely reducible $A_n(V)$-submodule of $X.$ Then $X/Z(n)$ is
an $A_{n-1}(V)$-module.

By the inductive assumption $X/Z(n)$ is a direct sum of
irreducible $A_{n-1}(V)$-modules. Without loss, we can assume that
$X/Z(n)$ is an irreducible $A_{n-1}(V)$-module. Then
 $X/Z(n)$ must be finite dimensional as $A_{n-1}(V)$ is semisimple.

On the other hand, $X^*$ is also an $A_n(V)$-module and
$(X/Z(n))^*$ is a submodule. Let $U=M(X)(0)^*.$ Then the Verma
type admissible module $M(U)=\bigoplus_{m\geq
0}A_{m,0}(V)\otimes_{A(V)}U$ is completely reducible and $M(U)(n)$
is a submodule of $X^*.$ As a result we see that $X^*=
M(U)(n)\oplus (X/Z(n))^*.$ Since $X$ is indecomposable, $X^*$ is
also indecomposable. This shows that $X/Z(n)=0$ and $X=Z(n)$ is
irreducible.

Here we give {\em another proof} for the statement that if $A(V)$ is semisimple
 and each irreducible admissible $V$-module is ordinary then $V$ is rational.
Again we will prove that $A_n(V)$ are semisimple for all $n.$ We
assume that $A_m(V)$ are semisimple for all $0\leq m< n.$

First by Lemma \ref{l2.3}, Lemma \ref{p2.4} and Corollary
\ref{co2.6} we have
\begin{eqnarray*}
& &A_{n,0}(V)\otimes_{A(V)} A_{0,n}(V)=\bigoplus_{i,j=1}^sM(U^i)(n)
\otimes (U^i)^*\otimes_{A(V)}U^j\otimes M((U^j)^*)(n)\\
& &\ \ \ =\bigoplus_{i=1}^sM(U^i)(n)\otimes M((U^i)^*)(n)\\
& &\ \ \ =\bigoplus_{i=1}^s\End(M(U^i)(n)).
\end{eqnarray*}
In particular, the dimension of $A_{n,0}(V)\otimes_{A(V)}
A_{0,n}(V)$ is $\sum_{i=1}^s(\dim M(U^i)(n))^2.$

By Theorem \ref{t2.8}, $A_{n,0}(V)*^n_{n,0}A_{0,n}(V)$ is a
homomorphic image of  $A_{n,0}(V)\otimes_{A(V)} A_{0,n}(V).$ So
its dimension is less than or equal to $\sum_{i=1}^s(\dim
M(U^i)(n))^2.$

From Lemma \ref{l2.8} we see that
\begin{eqnarray*}
& &\dim A_n(V)=\dim A_{n-1}(V)+\dim (O_n(V)/O_{n-1}(V)) \\
& &\ \ \ =\dim A_{n-1}(V)+\dim A_{n,0}(V)*^n_{n,0}A_{0,n}(V)\\
& &\ \ \ \leq \sum_{i=1}^s\sum_{m=0}^{n}(\dim M(U^i)(m))^2
\end{eqnarray*}
where we have used the fact that
$$\dim A_{n-1}(V)=\sum_{i=1}^s\sum_{m=0}^{n-1}(\dim M(U^i)(m))^2$$
 as $A_{n-1}(V)$ is semisimple.

On the other hand, $\{M(U^i)(m)|i=1,\cdots,s, m=0,\cdots,n\}$ are
the inequivalent irreducible modules for finite dimensional
algebra $A_n(V)$ by Theorem \ref{tha}. So the dimension of
$A_n(V)$ is at least $\sum_{i=1}^s\sum_{m=0}^{n}(\dim
M(U^i)(m))^2.$ This forces
$$\dim A_n(V)=\sum_{i=1}^s\sum_{m=0}^{n}(\dim M(U^i)(m))^2.$$
Thus
$$A_n(V)\cong \bigoplus_{i=1}^s\bigoplus_{m=0}^n\End(M(U^i)(m))$$
is semisimple. \qed

We remark that if the Verma type module $M(U)$ generated by any
irreducible $A(V)$-module $U$ is irreducible, $A(V)$ is semisimple
and $V$ is $C_2$-cofinite it is proved in  \cite{DLTYY} that  $V$
is rational. It is clear that the assumptions in Theorem
\ref{t2.11} are much weaker as any irreducible admissible module
is ordinary for any $C_2$-cofinite vertex operator algebra (cf.
\cite{KL}, \cite{Bu}, \cite{ABD}).

This theorem has several corollaries. Following \cite{KL} we call
$V$ $C_1$-cofinite if $V=\sum_{n\geq 0}V_n$ with $V_0=\C \1$
satisfies that $C_1(V)$ has finite codimension in $V$, where
$C_1(V)$ is spanned by $u_{-1}v$ and $L(-1)u$ for all $u,v\in
\sum_{n>0}V_n.$ It is proved in \cite{KL} that if $V$ is
$C_1$-cofinite then any irreducible admissible $V$-module is
ordinary.
\begin{coro} If $V$ is a $C_1$-cofinite simple vertex operator algebra
such that $A(V)$ is semisimple then $V$ is rational.
\end{coro}

A vertex operator algebra $V$ is called $C_2$-cofinite if the subspace
$C_2(V)$ spanned by $u_{-2}v$ for all $u,v\in V$ has finite codimension
in $V$ \cite{Z}. The $C_2$-cofiniteness has played a very important role
in the theory of vertex operator algebras and conformal field theory
(see \cite{Z}, \cite{DLM5}, \cite{GN}, \cite{NT}, \cite{DM2}, \cite{H},
\cite{DM3}).
It is proved in \cite{L} and \cite{ABD} that
$V$ is regular if and only if $V$ is rational and
$C_2$-cofinite. This implies the following corollary.

\begin{coro}\label{c2.10} Let $V$ be a simple $C_2$-cofinite vertex operator algebra.
Then the following are equivalent:

(a) $V$ is rational.

(b) $V$ is regular.

(c) $A(V)$ is semisimple.
\end{coro}

\section{An application: rationality of $V_{L}^{+}$}
\def\theequation{6.\arabic{equation}}
\setcounter{equation}{0}

We prove in this section the rationality of $V_L^+$ for all
positive definite even lattice $L.$ If the rank of $L$ is 1, the
rationality has been established previously in \cite{A1}.

Let $L$ be a positive definite even lattice of rank $d$ and
$V_{L}$ the vertex operator algebra associated to $L$ (cf.
\cite{B}, \cite{FLM}). Let $V_{L}^{+}$ be the fixed points of
$V_{L}$ under the automorphism $\theta$ lifted from the $-1$
isometry of $L$.  Then $V_{L}^{+}$ is a simple vertex operator
subalgebra of $V_{L}.$ If $d=1$ it is proved in \cite{A1} that
$V_{L}^{+}$ is rational. In this section we extend the rationality
to all $V_L^+.$ In \cite{DN2} and \cite{AD} we classify the
irreducible admissible modules of $V_{L}^{+}$. It turns out that
all the irreducible admissible modules are ordinary. Based on the
work \cite{DN2} and \cite{AD}, we prove in this section that
$A(V_{L}^{+})$ is semisimple. Thus by Theorem \ref{t2.11},
$V_{L}^{+}$ is rational.

We use the setting of \cite{FLM}. In particular, $\widehat{L}$ is
the canonical central extension of $L$ by the cyclic group
$<\kappa>$ of order 2:
$1\rightarrow<\kappa>\rightarrow\widehat{L}\rightarrow
L\rightarrow0$ with the commutator map
$c(\al,\be)=\kappa^{(\al,\be)}$ for $\al,\be\in L$. Let $e:
L\rightarrow\widehat{L}$ be a section such that $e_{0}=1$ and
$\epsilon: L\times L\rightarrow<\kappa>$  the corresponding
2-cocycle. We can assume that $\epsilon$ is bimultiplicative. Then
$\epsilon(\al,\be)\epsilon(\be,\al)=\kappa^{(\al,\be)}$,
$$
\epsilon(\al,\be)\epsilon(\al+\be,\gamma)=\epsilon(\be,\gamma)(\al,\be+\gamma),$$
and $e_{\al}e_{\be}=\epsilon(\al,\be)e_{\al+\be}$ for
$\al,\be,\gamma\in L$. Let $\theta$ denote the automorphism of
$\widehat{L}$ defined by $\theta(e_{\al})=e_{-\al}$ and
$\theta(\kappa)=\kappa$. Set
$K=\{a^{-1}\theta(a)|a\in\widehat{L}\}$.

Let $M(1)$ be the Heisenberg vertex operator algebra associated to
${\mathfrak h}=\C\otimes_{\Z}L.$ Then $V_L=M(1)\otimes \C[L]$
where $\C[L]$ is the group algebra of $L$ with a basis
$e^{\alpha}$ for $\alpha\in L$ and is an $\widehat L$-module such
that $e_\alpha e^\beta= \epsilon(\al,\be)e^{\alpha+\beta}.$

Recall that
$L^\circ=\{\,\lambda\in\h\,|\,(\alpha,\lambda)\in\Z\,\}$ is the
dual lattice of $L$.  There is an $\widehat{L}$-module structure
on $\C[L^\circ]=\bigoplus_{\lambda\in L^\circ}\C e^\lambda$ such
that $\kappa$ acts as $-1$ (see \cite{DL}). Let
$L^{\circ}=\cup_{i\in L^{\circ}/L}(L+\lambda_i)$ be the coset
decomposition such that $\lambda_0=0.$ Set
$\C[L+\lambda_i]=\bigoplus_{\alpha\in L}\C e^{\alpha+\lambda_i}.$
Then $\C[L^\circ]=\bigoplus_{i\in L^{\circ}/L}\C[L+\lambda_i]$ and
each $\C[L+\lambda_i]$ is an $\widehat L$-submodule of
$\C[L^\circ].$ The action of $\widehat L$ on $\C[L+\lambda_i]$ is
as follows:
$$e_{\alpha}e^{\beta+\lambda_i}=\e(\a,\b)e^{\a+\b+\l_i}$$
for $\alpha,\,\b\in L.$ On the surface, the module structure on
each $\C[L+\lambda_i]$ depends on the choice of $\lambda_i$ in
$L+\lambda_i.$ It is easy to prove that different choices of
$\lambda_i$ give isomorphic $\widehat L$-modules.

 Set $\C[M]=\bigoplus_{\lambda\in M}\C e^{\lambda}$ for a subset $M$ of
$L^{\circ}$, and define $V_M=M(1)\otimes\C[M]$. Then $V_L$ is a
rational vertex operator algebra and $V_{L+\lambda_i}$ for $i\in
L^{\circ}/L$ are the irreducible modules for $V_L$ (see \cite{B},
\cite{FLM}, \cite{D1}, \cite{DLM1}).

Define a linear isomorphism
$\theta:V_{L+\lambda_i}\to V_{L-\lambda_i}$ for $i\in L^{\circ}/L$ by
\begin{align*}
\theta(\beta_{1}(-n_{1})\beta_{2}(-n_{2})\cdots \beta_{k}(-n_{k})
e^{\alpha+\lambda_i})=(-1)^{k}\beta_{1}(-n_{1})\beta_{2}(-n_{2})\cdots \beta_{k}(-n_{k})
e^{-\alpha-\lambda_i}
\end{align*}
for $\beta_i\in\h,\,n_i\geq1$ and $\alpha\in L$
if $2\lambda_i\not\in L,$ and
\begin{multline*}
\theta(\beta_{1}(-n_{1})\beta_{2}(-n_{2})\cdots \beta_{k}(-n_{k})e^{\alpha+\lambda_i})\\
=(-1)^{k}c_{2\l_i}\e(\alpha,2\lambda_i)\beta_{1}(-n_{1})\beta_{2}(-n_{2})\cdots \beta_{k}(-n_{k})e^{-\alpha-\lambda_i}
\end{multline*}
if $2\lambda_i\in L$ where $c_{2\l_i}$ is a square root of
$\e(2\l_i,2\l_i).$
  Then $\theta$ defines a linear
isomorphism from $V_{L^{\circ}}$ to itself such that
 $$\theta Y(u,z)v=Y(\theta u,z)\theta v$$
for $u\in V_L$ and $v\in V_{L^{\circ}}.$
In particular, $\theta$ is an automorphism of $V_L$ which induces
an automorphism of $M(1).$

For any $\theta$-stable subspace $U$ of $V_{L^{\circ}}$, let
$U^\pm$ be the $\pm1$-eigenspace of $U$ for $\theta$. Then $V_L^+$
is a simple vertex operator algebra.

Also recall the $\theta$-twisted Heisenberg algebra $\h[-1]$ and
its irreducible module $M(1)(\theta)$ from \cite{FLM} and
\cite{AD}.  Let $\chi$ be a central character of $\widehat{L}/K$
such that $\chi(\kappa)=-1$ and $T_{\chi}$ the irreducible
$\widehat{L}/K$-module with central character $\chi$. Then
$V_L^{T_{\chi}}$ is an irreducible $\theta$-twisted $V_L$-module
(see \cite{FLM}, \cite{D2} and \cite{DL}). We also define an
action of $\theta$ on  $V_L^{T_{\chi}}$ such that \begin{align*}
\theta(\beta_{1}(-n_{1})\beta_{2}(-n_{2})\cdots
\beta_{k}(-n_{k})t)=(-1)^{k}\beta_{1}(-n_{1})\beta_{2}(-n_{2})\cdots
\beta_{k}(-n_{k})t
\end{align*}
for $\beta_i\in\h,\,n_i\in \frac{1}{2}+\Z_{+}$ and $t\in
T_{\chi}$. Recall that $L^{\circ}=\cup_{i\in
L^{\circ}/L}(L+\lambda_i).$ Here is the classification of
irreducible modules
 for $V_L^+$ (see
\cite{DN2} and \cite{AD}).
\begin{theorem}
Let $L$ be a positive definite even lattice. Then any irreducible
admissible $V_L^+$-module is isomorphic to one of irreducible
modules $V_L^{\pm},\, V_{L+\l_i}(2\lambda_i\notin L),$
$V_{L+\lambda_i}^{\pm}\,(2\lambda_i\in L)$
 and $V_L^{T_{\chi},\pm}$ for any irreducible $\widehat{L}/K$-module $T_{\chi}$ with central character $\chi$.
\end{theorem}

Then the irreducible $A(V_L^+)$-modules are the top levels $W(0)$
of irreducible $V_L^+$-modules $W$ given as follows:
\begin{align*}
&\charge{+}(0)=\C\1,\quad
\charge{-}(0)=\h(-1)\bigoplus(\bigoplus_{\alpha\in L_2}\C
(e^{\alpha}-e^{-\alpha})),\\
&\charlam{\lambda_i}(0)=\bigoplus_{\alpha\in\Delta(\lambda_i)}\C e^{\lambda_i+\alpha}\quad(2\lambda_i\notin L),\\
&\charlam{\lambda_i}^{\pm}(0)=\sum_{\alpha\in\Delta(\lambda_i)}\C(e^{\lambda_i+\alpha}\pm \theta e^{\lambda_i+\alpha})\,
\quad(2\lambda_i\in L),\\
&\charge{T_{\chi},+}(0)=T_{\chi},\quad\charge{T_{\chi},-}(0)=\h(-1/2)\otimes T_{\chi}.
\end{align*}
Here $\h(-1)=\{h(-1)|h\in\h\}\subset M(1)$ and
$\h(-1/2)=\{h(-1/2)|h\in\h\}\subset M(1)(\theta).$

Let $\{h_1,\cdots,h_d\}$ be an orthonormal basis of $\h.$ Recall
from \cite{DN3} and \cite{AD} the following vectors in $V_L^+$ for
$a,b=1,\cdots,d$ and $\alpha\in L$
\begin{align*}
S_{ab}&(m,n)=h_a(-m)h_b(-n),\\
E^u_{ab}&=5 S_{ab}(1,2)+25 S_{ab}(1,3)+36 S_{ab}(1,4)+16 S_{ab}(1,5)\,(a\neq b),\\
\bar{E}^u_{ba}&=S_{ab}(1,1)+14S_{ab}(1,2)+41S_{ab}(1,3)+44S_{ab}(1,4)+16 S_{ab}(1,5)\,(a\neq b),\\
E_{aa}^u&=E^{u}_{ab}E^{u}_{ba},\\
E^t_{ab}&=-16(3 S_{ab}(1,2)+14S_{ab}(1,3)+19S_{ab}(1,4)+8 S_{ab}(1,5))\,(a\neq b),\\
\bar{E}^t_{ba}&=-16(5S_{ab}(1,2)+18 S_{ab}(1,3)+21 S_{ab}(1,4)+8 S_{ab}(1,5))\,(a\neq b),\\
E_{aa}^t&=E^{t}_{ab}E^{t}_{ba},\\
\Lambda_{ab}&=45 S_{ab}(1,2)+190 S_{ab}(1,3)+240 S_{ab}(1,4)+96 S_{ab}(1,5),\\
E^{\a}&=e^{\a}+e^{-\a}.
\end{align*}
For $v\in V_L^+$ we denote $v+O(V_L^+)$ by $[v].$ Let $A^u$ and
$A^t$ be the linear subspace of $A(M(1)^+)$ spanned by
$E_{ab}^u+O(M(1)^+)$ and $E_{ab}^{t}+O(M(1)^+)$ respectively for
$1\leq a,\,b\leq d.$ Then $A^t$ and $A^u$ are two sided ideals of
$A(M(1)^+).$ Note that the natural algebra homomorphism from
$A(M(1)^+)$ to $A(V_L^+)$ gives embedding of $A^u$ and $A^t$ into
$A(V_L^+).$ We should remark that the $A^u$ and $A^t$ are
independent of the choice of the orthonormal basis
$\{h_1,\cdots,h_d\}.$

By Lemma 7.3 of \cite{AD} we know that
$$V_L^-(0)=\h(-1)\bigoplus(\sum_{\a\in L_2}\C[E^\a]\a(-1)),$$
where $L_2=\{\alpha\in L|(\a,\a)=2\}.$ Let
$L_{2}=\{\pm\al_{1},\cdots,\pm\al_{r},\pm\al_{r+1},\cdots,\pm\al_{r+l}\}$
be such that $\{\al_{1},\cdots,\al_{r}\}$ are linearly independent
and $\{\al_{r+1},\cdots,\al_{r+l}\} \ \subseteq
\bigoplus_{i=1}^{r}\Z_{+}\al_{i}.$ We can choose the orthonormal
basis $\{h_i | \ i=1,\cdots,d\}$ so that $h_{i}\in
\C\al_{1}+\cdots+\C\al_{i}$, for $i=1,\cdots,r.$ Then we have
$$\al_{i}(-1)=a_{i1}h_{1}(-1)+\cdots+a_{ii}h_{i}(-1), \
i=1,\cdots,r,$$
$$
\al_{j}(-1)=a_{j1}h_{1}(-1)+\cdots +a_{jr}h_{r}(-1), \
j=r+1,\cdots,r+l,$$ where $a_{ii}\neq 0, i=1,\cdots,r$. For
$i\in\{1,2,\cdots,l\}$, let $k_{i}$ be such that
$$a_{r+i,k_{i}}\neq 0, \ a_{r+i,k_{i}+1}=\cdots=a_{r+i,r}=0.$$

We know from \cite{AD} that $e^i=h_i(-1)$ for $i=1,\cdots,d$ and
$e^{d+j}=[E^{\a_j}]\alpha_j(-1)$ for $j=1,\cdots,r+l$ form a basis
of $V_L^-(0).$ We first construct a two-sided ideal of $A(V_L^+)$
isomorphic to $\End(V_L^-(0)).$ Recall $E^u_{ij}$ for
$i,j=1,\cdots,d.$ We now extend the definition of $E^u_{ij}$  to
all $i,j=1,\cdots,d+r+l$  and the linear span of $E^u_{ij}$ will
be the ideal of $A(V_L^+)$ isomorphic to $\End V_L^-(0)$ (with
respect to the basis $\{e^1,\cdots,e^{d+r+l}\}$ ).

For the notational convenience,  we also write $E^u_{i,j}$ for
$E^u_{ij}$ from now on. Define
$$
[E_{j,d+i}^{u}]=\frac{1}{4\epsilon(\al_{i},\al_{i})a_{ii}}[E_{ji}^{u}*E^{\al_{i}}],
\ i=1,\cdots,r,j=1,\cdots,d,$$
$$
[E_{j,d+r+i}^{u}]=\frac{1}{4\epsilon(\al_{r+i},\al_{r+i})a_{r+i,k_{i}}}[E_{j,k_{i}}^{u}*E^{\al_{r+i}}],
\ i=1,\cdots,l,j=1,\cdots,d.$$  Define
$$
[E_{d+i,j}^{u}]=\sum\limits_{k=1}^{r}a_{ik}[E^{\al_{i}}]*[E_{kj}^{u}],
\ i=1,\cdots,r+l,j=1,\cdots,d,$$ where  $a_{ij}=0$, for $1\leq
i<j\leq r$.
 Recall from \cite{DN3} and
\cite{AD} that $[E^u_{ab}]h_c(-1)=\delta_{c,b}h_a(-1)$ for
$a,b,c=1,\cdots,d.$
\begin{lem}\label{l5.1} The following holds:
$$[E^u_{ij}]e^k=\delta_{k,j}e^i,\ [E^u_{st}]e^k=\delta_{t,k}e^s$$
for $i,t=1,\cdots,d,$ $j,s=d+1,\cdots,d+r+l$ and
$k=1,\cdots,d+r+l.$
\end{lem}

\pf Let $h\in {\mathfrak h}$ such that $(h,h)\ne 0.$ Then $\omega_h=
\frac{1}{2(h,h)}h(-1)^2$ is a Virasoro element with central charge
1. Note that $\omega_h\beta(-1)=\frac{(\beta,h)^4}{2(h,h)}h(-1)$
for any $\beta\in {\mathfrak h}.$ For $\alpha\in L_{2}$ then
$[E^\alpha]*[E^\alpha]=4\e(\a,\a)[\omega_\alpha]$ in $A(V_L^+)$ by
Proposition 4.9 of \cite{AD}.  Then for $i=1,\cdots,d,
j=1,\cdots,r$, we have
\begin{eqnarray*}
& & [E_{i,d+j}^{u}]e^{d+j}=[E_{i,d+j}^{u}]([E^{\al_{j}}]\al_{j}(-1))\\
&
&\ \ \ =\frac{1}{4\epsilon(\al_{j},\al_{j})a_{jj}}([E_{ij}^{u}]*[E^{\al_{j}}]*[E^{\al_{j}}])\al_{j}(-1)\\
& &\ \ \ =\frac{1}{a_{jj}}[E_{ij}^{u}]\al_{j}(-1)=h_{i}(-1).
\end{eqnarray*}

Let $k\in\{1,\cdots,r+l\}$ such that $k\ne j.$ Then
\begin{eqnarray*}
& &[E_{i,d+j}^{u}]e^{d+k}=[E_{i,d+j}^{u}]([E^{\al_{k}}]\al_{k}(-1))\\
&
&\ \ \ =\frac{1}{4\epsilon(\al_{j},\al_{j})a_{jj}}([E_{ij}^{u}]*[E^{\al_{j}}]*[E^{\al_{k}}])\al_{k}(-1).
\end{eqnarray*}
By Proposition 5.4 of \cite{AD}, we have
$$[E^{\al_{j}}]*[E^{\al_{k}}]=\sum\limits_{p}[v^{p}]*[E^{\al_{j}+\al_{k}}]*[w^{p}]+\sum\limits_{q}[x^{q}]*
[E^{\al_{j}-\al_{k}}]*[y^{q}],
$$
where $v^{p}, w^{p}, x^{q}, y^{q}\in M(1)^{+}$.
 Since $A^{u}$ is an ideal of
$A(M(1)^{+})$, we have $[E_{ji}^{u}]*[v^p], [E_{ji}^{u}]*[x^q]\in
A^{u}.$ By the proof of Proposition 7.2 of \cite{AD}, we know that
$A^{u}[E^{\al}]\al(-1)=0$, for any $\al\in L_{2}$. So by Lemma
7.1 and Proposition 7.2 of \cite{AD}, we have
$$[E_{i,d+j}^{u}]e^{d+k}=0, \ i=1,\cdots,d,
j=1,\cdots,r, k=1,\cdots,r+l, j\neq k.$$ It follows from the proof
of Proposition 7.2 of \cite{AD} that
$$E^u_{i,j+d}e^s=[E_{ij}^u]*[E^{\a_j}]h_s(-1)=0$$
for $s=1,\cdots,d$ as $[E^{\a_j}]h_s(-1)\in\sum_{p=1}^{r+l}\C
[e^{\a_{p}}-e^{-\a_{p}}].$ This completes the proof for
$E^u_{i,j+d}$ for $i=1,\cdots,d,$ $j=1,\cdots,r.$ The other cases
can be done similarly. \qed

Recall $H_a$ and $\omega_a=\omega_{h_a}$ for $a=1,\cdots,d$
 from \cite{DN3} and \cite{AD}. The following lemma collects some formulas
from Propositions 4.5, 4.6, 4.8 and 4.9 of \cite{AD}.
\begin{lem}\label{adl} For any indices $a,\,b,\,c,\,d$,
\begin{align}
&[\w_a]*[E^{u}_{bc}]=\delta_{ab}[E^{u}_{bc}],\label{e-1}\\
&[E^{u}_{bc}]*[\w_a]=\delta_{ac}[E^{u}_{bc}].\label{e-2}\\
&[E_{ab}^{u}]*[E_{cd}^{t}]=[E_{cd}^{t}]*[E_{ab}^{u}]=0,
\label{equ1}\\
&[\La_{ab}]*[E_{cd}^{u}]=[\La_{ab}]*[E_{cd}^{t}]=[E_{cd}^{u}]*[\La_{ab}]=[E_{cd}^{t}]*[\La_{ab}]=0
\ (a\neq b),\label{equ2}
\end{align}
For distinct $a,b,c,$
\begin{align}
&\left(70 [H_a]+1188[\w_a]^2-585 [\w_a]+27\right)*[H_a]=0,\label{equation1}\\
&([\w_a]-1)*\left([\w_a]-\frac{1}{16}\right)*\left([\w_a]-\frac{9}{16}\right)*[H_a]=0,\label{equation2}\\
&-\frac{2}{9}[H_a]+\frac{2}{9}[H_b]=2[E_{aa}^u]-2[E_{bb}^u]+\frac{1}{4}[E_{aa}^t]-\frac{1}{4}[E_{bb}^t],
\label{equation3}\\
\begin{split}\label{equation4}
&-\frac{4}{135}(2[\w_a]+13)*[H_a]+\frac{4}{135}(2[\w_b]+13)*[H_b]\\
&\ \ \
=4([E_{aa}^u]-[E_{bb}^u])+\frac{15}{32}([E_{aa}^t]-[E_{bb}^t]),
\end{split}\\
&[\w_b]*[H_a]=-\frac{2}{15}([\w_a]-1)*[H_a]+\frac{1}{15}([\w_b]-1)*[H_b],\label{equation5}\\
&[\Lambda_{ab}]^2=4[\w_a]*[\w_b]-\frac{1}{9}([H_a]+[H_b])-([E_{aa}^u]+[E_{bb}^u])-\frac{1}{4}
([E_{aa}^t]+[E_{bb}^t]),\label{equation6}\\
&[\Lambda_{ab}]*[\Lambda_{bc}]=2[\w_b]*[\Lambda_{ac}].\label{equation7}
\end{align}
For $\alpha\in L$ such that $(\alpha,\alpha)=2k\ne 2,$
\begin{align}
&[H_{\alpha}]*[E^{\alpha}]=\frac{18(8k-3)}{(4k-1)(4k-9)}\left([\w_{\alpha}]-\frac{k}{4}\right)\left([\w_{\alpha}]-
\frac{3(k-1)}{4(8k-3)}\right)[E^{\alpha}],\label{equation8}\\
&\left([\w_{\alpha}]-\frac{k}{4}\right)\left([\w_{\alpha}]-\frac{1}{16}\right)\left([\w_{\alpha}]-
\frac{9}{16}\right)[E^{\alpha}]=0.\label{equation9}
\end{align}
If $\alpha\in L_2,$
\begin{align}
&[E^\alpha]*[E^\alpha]=4\e(\a,\a)[\w_\alpha],\label{equation10}\\
&[H_\alpha]*[E^\alpha]+[E^\alpha]*[H_\alpha]=-12[\w_\alpha]*\left([\w_{\alpha}]-\frac{1}{4}\right)*[E^{\alpha}],
\label{equation11}\\
&([\w_{\alpha}]-1)*\left([\w_{\alpha}]-\frac{1}{4}\right)*\left([\w_{\alpha}]-\frac{1}{16}\right)*\left([\w_{\alpha}]-
\frac{9}{16}\right)*[E^{\alpha}]=0.\label{equation12}
\end{align}
For any $\al\in L$,
\begin{align}
I^{t}*[E^{\al}]=[E^{\al}]*I^{t},\label{equation13}
\end{align}
where $I^{t}$ is the identity of the simple algebra $A^{t}$.
\end{lem}

\begin{lem}\label{l5.2} For any $\al\in L_{2}$, we have
$$A^u*[E^{\al}]*A^u=0.
$$
\end{lem}

\pf Let $\al\in L_{2}$ and $\{h_{1},\cdots,h_{d}\}$ be  an
orthonormal basis of ${\f h}$ such that $h_{1}\in\C\al$. ($A^u$ is
independent of the choice of orthonormal basis.) By
(\ref{e-1})-(\ref{e-2}) and (\ref{equation12}), we have
$[E^{\alpha}]=f([\omega_{\a}])*[E^{\alpha}]=[E^{\alpha}]*f([\omega_{\a}])$
for some polynomial $f(x)$ with $f(0)=0.$ Note that
$\omega_{\alpha}=\omega_1.$ By (\ref{e-1})-(\ref{e-2}), we  only
need to prove that
$$
[E_{i1}^{u}]*[E^{\al}]*[E_{1s}^{u}]=0, \ i,s=1,2,\cdots,d.$$ Let
$a=1$, $b\neq 1$ in (\ref{equation5}). Multiplying
(\ref{equation5}) by $[E_{i1}^{u}]$ on left and using (\ref{e-2})
and (\ref{equ1}), we have
$$[E_{i1}^{u}]*[H_{b}]=0, \ b\neq 1.$$
Then setting $a=1$, $b\neq 1$ in (\ref{equation3}) and multiplying
(\ref{equation3}) by $[E_{i1}^{u}]$ on left yields
$$[E_{i1}^{u}]*[H_{1}]=-9[E_{i1}^{u}].$$
Let $a=1$, $b\neq 1$ in (\ref{equation6}). Multiplying
(\ref{equation6}) by $[E_{1s}^{u}]$ on right and using (\ref{e-1})
and (\ref{equ2}), we have
$$
-\frac{1}{9}[H_{1}]*[E_{1s}^{u}]-\frac{1}{9}[H_{b}]*[E_{1s}^{u}]=[E_{1s}^{u}].$$
On the other hand, multiplying (\ref{equation3}) by $[E_{1s}^{u}]$
on right yields
$$
-\frac{1}{9}[H_{1}]*[E_{1s}^{u}]+\frac{1}{9}[H_{b}]*[E_{1s}^{u}]=[E_{1s}^{u}].$$
Comparing the above two formulas, we have
$$[H_{1}]*[E_{1s}^{u}]=-9[E_{1s}^{u}], \
[H_{a}]*[E_{1s}^{u}]=0, \ a\neq 1.$$ So
$$
[E_{i1}^{u}]*[H_{1}]*[E^{\al}]*[E_{1s}^{u}]+[E_{i1}^{u}]*[E^{\al}]*[H_{1}]*
[E_{1s}^{u}]=-18[E_{i1}^{u}]*[E^{\al}]*[E_{1s}^{u}].$$
But by (\ref{e-1})-(\ref{e-2}) and (\ref{equation11})
 we have
$$
[E_{i1}^{u}]*[H_{1}]*[E^{\al}]*[E_{1s}^{u}]+[E_{i1}^{u}]*[E^{\al}]*[H_{1}]*
[E_{1s}^{u}]=-9[E_{i1}^{u}]*[E^{\al}]*[E_{1s}^{u}].$$ This implies
that $ [E_{i1}^{u}]*[E^{\al}]*[E_{1s}^{u}]=0,$ as required.
 \qed

We now define $E^u_{i,j}$ for all $i,j=1,\cdots,d+r+l.$  Set
$$
[E_{d+i,d+j}^{u}]=[E_{d+i,1}^{u}]*[E_{1,d+j}^{u}], \
i,j=1,\cdots,r+l.$$ It is easy to see that
$[E_{d+i,1}^{u}]*[E_{1,d+j}^{u}]=[E_{d+i,k}^{u}]*[E_{k,d+j}^{u}]$,
$k=2,\cdots,d$.

Denote by $A_{L}^u$ the subalgebra of $A(V_{L}^{+})$ generated by
$\{[E_{ij}^{u}],[E_{d+p,j}^{u}],[E_{i,d+p}^{u}]|i,j=1,\cdots,d, \
p=1,\cdots,r+l\}$. From Lemma \ref{l5.2}, (\ref{equation10}) and
the definition of $[E_{ij}^{u}]$, $i,j=1,\cdots,d+r+l$, we can
easily deduce the following result.

\begin{lem}\label{l5.4} $A_{L}^u$ is a matrix algebra over $\C$ with basis  $\{[E_{ij}^{u}]|i,j=1,\cdots,d+r+l\}$
such that
$$
[E_{ij}^{u}]*[E_{ks}^{u}]=\de_{j,k}[E_{is}^{u}], \
[E_{ij}^{u}]e^{k}=\de_{j,k}e^{i}, \ i,j,k,s=1,2,\cdots,d+r+l.$$
\end{lem}

\begin{lem}\label{l5.3} Let $\al\in L$ and $\al\notin L_{2}$, then
$$[E^{\al}]*A^u=0.$$
\end{lem}

\pf Let $\{h_{1},\cdots,h_{d}\}$ be  an orthonormal basis of ${\f
h}$ such that $h_{1}\in\C\al$. If $|\al|^{2}=2k$ and $k\neq 4,$
the lemma follows from
(\ref{e-1})-(\ref{e-2}) and (\ref{equation9}).

If $|\al|^{2}=8,$ by (\ref{e-1})-(\ref{e-2}) and (\ref{equation9})
we have
$$
[E_{ab}^{u}]*[E^{\al}]=[E^{\al}]*[E_{ba}^{u}]=0, $$ for all $1\leq
a,b\leq d$ and $b\neq 1$. By (\ref{e-1}) and (\ref{equation8}), we
have
\begin{equation}\label{ea5.20}
[E_{a1}^{u}]*[H_{1}]*[E^{\al}]=0.
\end{equation}
On the other hand, for $a\neq 1$, by
(\ref{equation3})-(\ref{equation4}) and (\ref{equ1}), we have
$$
-\frac{2}{9}[E_{a1}^{u}]*([H_{a}]-[H_{1}])*[E^{\al}]=-2[E_{a1}^{u}]*[E^{\al}],$$
$$
\frac{4}{135}[E_{a1}^{u}]*(-13[H_{a}]+15[H_{1}])*[E^{\al}]=-4[E_{a1}^{u}]*[E^{\al}].$$
Therefore by (\ref{ea5.20}), we have
$$\frac{1}{9}[E_{a1}^{u}]*[H_{a}]*[E^{\al}]=[E_{a1}^{u}]*[E^{\al}], \ a\neq 1,$$
$$\frac{13}{135}[E_{a1}^{u}]*[H_{a}]*[E^{\al}]=[E_{a1}^{u}]*[E^{\al}], \ a\neq 1.$$
This means that
$$
[E_{a1}^{u}]*[E^{\al}]=0.$$ Since
$[H_{\al}]*[E^{\al}]=[E^{\al}]*[H_{\al}]$, we similarly have
$$[E^{\al}]*[E_{1a}^{u}]=0.$$
This completes the proof. \qed

\begin{lem}\label{l5.5} $A_{L}^{u}$ is an ideal of $A(V_{L}^{+})$.
\end{lem}

\pf By Proposition 5.4  of \cite{AD}, (\ref{equ1}),
(\ref{equation13}) and Lemmas \ref{l5.2}-
 \ref{l5.3}, it is enough to prove that $[E^{\al_{i}}]*[E_{jk}^{u}], [E_{jk}^{u}]*[E^{\al_{i}}]\in A_{L}^u$,
  $j,k=1,\cdots,d, \
 i=1,\cdots,r+l$.

 Let $\al\in L_{2}$. For convenience, let $a_{ij}=0$, for
$1\leq i<j\leq r$ and $k_{i}=i$ for $1\leq i\leq r$. Since
$\al_{i}(-1)=\sum\limits_{k=1}^{d}a_{ik}h_{k}(-1)$, we have
$$\omega_{\al_{i}}=\frac{1}{4}\sum\limits_{k=1}^{d}a_{ik}^{2}h_{k}(-1)^{2}+\frac{1}{4}\sum\limits_{p\neq
q}a_{ip}a_{iq}h_{p}(-1)h_{q}(-1).$$ Recall from \cite{AD} that
$$
[S_{ab}(1,1)]=[E_{ab}^{u}]+[E_{ba}^{u}]+[\La_{ab}]+\frac{1}{2}[E_{ab}^{t}]+\frac{1}{2}[E_{ba}^{t}],
\ a\neq b.$$ So from (\ref{equ1}) and (\ref{equ2}) we have
$$
[E_{jk}^{u}]*[\omega_{\al_{i}}]=\frac{1}{2}a_{ik}^{2}[E_{jk}^{u}]+\frac{1}{2}\sum\limits_{p\neq
k}a_{ik}a_{ip}[E_{jp}^{u}], \ j,k=1,\cdots,d, \ i=1,\cdots,r+l.$$
Then it can  easily be deduced that
$$(a_{ik}[E_{jk_{i}}^{u}]-a_{ik_{i}}[E_{jk}^{u}])*[\omega_{\al_{i}}]=0,
\ j,k=1,\cdots,d, \ i=1,\cdots,r+l.$$ Then by (\ref{equation12}),
we have
$$(a_{ik}[E_{jk_{i}}^{u}]-a_{ik_{i}}[E_{jk}^{u}])*[E^{\al_{i}}]=0,
\ j,k=1,\cdots,d, \ i=1,\cdots,r+l.$$ So
$[E_{jk}^{u}]*[E^{\al_{i}}]\in A_{L}^{u}$, for all
$j,k=1,\cdots,d, i=1,\cdots,r+l.$ Similarly, we have
$$[E^{\al_{i}}]*[E_{jk}^{u}]=0, \ k=1,\cdots,d, \ i=1,\cdots,r+l, \ j=r+1,\cdots,d,$$
$$
[E^{\al_{i}}]*[\sum\limits_{b=1}^{r}a_{kb}E_{bj}^{u}]=\frac{(\al_{i},\al_{k})}{2}[E^{\al_{i}}]*[\sum\limits_{b=1}^{r}
a_{ib}E_{bj}^{u}], \   j=1,\cdots,d, \ k=1,\cdots,r, \
i=1,\cdots,r+l.$$ Since
 both $\{\al_{1},\cdots,\al_{r}\}$ and $\{h_{1},\cdots,h_{r}\}$
 are linearly independent, it follows that for each
 $i=1,\cdots,r, \ j=1,\cdots,d$,
 $[E_{ij}^{u}]$ is a linear combination of $a_{11}[E_{1j}^{u}]$,
 $[a_{21}E_{1j}+a_{22}E_{2j}^{u}], \cdots,
 [a_{r1}E_{1j}^{u}+\cdots+a_{rr}E_{rj}^{u}]$. Therefore
 $[E^{\al_{i}}]*[E_{jk}^{u}]\in A_{L}^u$, $j,k=1,\cdots,d, \
 i=1,\cdots,r+l$. \qed

 For $0\neq \al\in L$, let $\{h_{1},\cdots,h_{d}\}$ be  an
orthonormal basis of ${\f h}$ such that $h_{1}\in\C\al$. Define
$$[B_{\al}]=2^{|\al|^{2}-1}([I^{t}]*[E^{\al}]-\frac{2|\al|^{2}}{2|\al|^{2}-1}[E_{11}^{t}]*[E^{\al}]),$$
and $[B_{0}]=[I^{t}]$
(see formula (6.5) of \cite{AD}).
\begin{lem}\label{l5.6}
For $\al\in L$, $[E_{ij}^{t}]\in A^{t}$,
$[B_{\al}]*[E_{ij}^{t}]=[E_{ij}^{t}]*[B_{\al}]$.
\end{lem}

\pf It is enough to prove that
$$
[B_{\al}]*[E_{ij}^{t}]=[E_{ij}^{t}]*[B_{\al}],$$ for $i=1$ or
$j=1$. By the definition of $[E_{ab}^{t}]$ and the fact that
$$[I^{t}]*[E^{\al}]=[E^{\al}]*[I^{t}]$$
and
$$
[I^{t}]*[\La_{ab}]=[I^{t}]*[E_{ab}^{u}]=0, \ a\neq b,$$ we have
$$
[B_{\al}]*[E_{ab}^{t}]=[B_{\al}]*(-[S_{ab}(1,1)]-2[S_{ab}(1,2)]),$$
$$
[E_{ab}^{t}]*[B_{\al}]=(-[S_{ab}(1,1)]-2[S_{ab}(1,2)])*[B_{\al}].$$
 Let $b\neq 1$. Similar to
the proof of Lemma 7.5 of \cite{AD}, we have
\begin{eqnarray*}
& & \ \ \
(2|\al|^{2}-1)([E_{1b}^{t}]+3[E_{1b}^{u}]+[\La_{1b}])*[E^{\al}]+[E^{\al}]*([E_{1b}^{t}]+3[E_{1b}^{u}]+[\La_{1b}])\\
& &
=-([E_{b1}^{t}]-[E_{1b}^{u}]+[\La_{1b}])*[E^{\al}]-(2|\al|^{2}-1)[E^{\al}]*([E_{b1}^{t}]-[E_{1b}^{u}]+[\La_{1b}])
\end{eqnarray*}
\begin{eqnarray*}
& & \ \ \
(2|\al|^{2}-1)(\frac{1}{16}[E_{1b}^{t}]*[E^{\al}]+[\omega_{b}]*[\La_{1b}]*[E^{\al}])
+\frac{1}{16}[E^{\al}]*[E_{1b}^{t}]+[E^{\al}]*[\omega_{b}]*[\La_{1b}]\\
& &
=\!-\!\frac{9}{16}[E_{b1}^{t}]*[E^{\al}]\!+\![\omega_{b}]*[\La_{1b}]*[E^{\al}]
\!-\!(2|\al|^{2}-1)(\frac{9}{16}[E^{\al}]*[E_{b1}^{t}]+[E^{\al}]*[\omega_{b}]*[\La_{1b}]).
\end{eqnarray*}
So we have
$$(2|\alpha|^{2}-1)[E_{1b}^{t}]*[E^{\al}]=-[E^{\al}]*[E_{1b}^{t}]+x,$$
where $x\in
A_{L}^{u}+\C[\La_{1b}]*[E^{\al}]+\C[E^{\al}]*[\La_{1b}]+\C[\omega_{b}]*[\La_{1b}]*[E^{\al}]
+\C[E^{\al}]*[\La_{1b}]*[\omega_{b}]$. Since $y*x=0$ for any $y\in
A^{t}$, we have
\begin{eqnarray*}
& &  \ \ \
[B_{\al}]*[E_{1b}^{t}]\\
 &
 &=2^{|\al|^{2}-1}\left(-(2|\al|^{2}-1)[E_{1b}^{t}]*[E^{\al}]+2|\al|^{2}[E_{11}^{t}]*[E_{1b}^{t}]*[E^{\al}]\right)\\
 & &=2^{|\al|^{2}-1}[E_{1b}^{t}]*[E^{\al}]=[E_{1b}^{t}]*[B_{\al}].
\end{eqnarray*}
Similarly,
$$[B_{\al}]*[E_{b1}^{t}]=[E_{b1}^{t}]*[B_{\al}],$$
completing the proof. \qed

\begin{lem}\label{l5.7} $A_{L}^{t}$ is an ideal of $A(V_{L}^{+})$ and
$A_{L}^{t}\cong A^{t}\otimes_{\C} {\C}[\widehat{L}/ K]/J$, where
${\C}[\widehat{L}/K]$ is the group algebra of $\widehat{L}/ K$ and
$J$ is the ideal of ${\C}[\widehat{L}/ K]$ generated by $\kappa
K+1.$
\end{lem}

\pf By Proposition 5.4 of \cite{AD} and Lemmas
\ref{l5.5}-\ref{l5.6}, it is easy to check that $A_{L}^{t}$ is an
ideal of $A(V_{L}^{+})$. Similar to the proof of Proposition 7.6
of \cite{AD}, we have
$$
[B_{\al}]*[B_{\be}]=\epsilon(\al,\be)[B_{\al+\be}],$$ for
$\al,\be\in L$ where $\epsilon(\al,\be)$ is understood to be
$\pm 1$ by identifying $\kappa$ with $-1.$
Then the lemma follows from Proposition 7.6 of
\cite{AD} and Lemma \ref{l5.6}. \qed

It is clear that $A_{L}^{u}\cap A_{L}^{t}=0$. Let
$$
\bar{A}(V_{L}^{+})=A(V_{L}^{+})/(A_{L}^{u}\oplus A_{L}^{t}),
$$
and for $x\in {A}(V_{L}^{+})$, we still denote the image of $x$ in
$\bar{A}(V_{L}^{+})$ by $x$.

\begin{lem}\label{l5.8}
In $\bar{A}(V_{L}^{+})$, we have
\begin{equation}\label{ea5.22}
[H_{a}]=[H_{b}], \ 1\leq a,b\leq d,
\end{equation}
\begin{equation}\label{ea5.23}
([\omega_{a}]-\frac{1}{16})*[H_{a}]=0, \ 1\leq a\leq d,
\end{equation}
\begin{equation}\label{ea5.24}
\frac{128}{9}[H_{a}]*\frac{128}{9}[H_{a}]=\frac{128}{9}[H_{a}], \
1\leq a\leq d,
\end{equation}
\begin{equation}\label{ea5.25}
[\La_{ab}]*[H_{c}]=0, \ 1\leq a,b,c\leq d, \ a\neq b.
\end{equation}
\end{lem}

\pf (\ref{ea5.22}) follows from (\ref{equation3})  and
(\ref{ea5.23}) follows from (\ref{equation4}) and
(\ref{equation5}). Then from (\ref{equation1}) we can
get (\ref{ea5.24}). By (\ref{equation6}), we have
$$
[\La_{ab}]^{2}*[H_{c}]=0, \ a\neq b.
$$
If $d\geq 3$, then by (\ref{ea5.22}) we can let $c\neq a$, $c\neq
b$. So by (\ref{equation7}) and (\ref{ea5.23}),
\begin{eqnarray*}
& & \   \ \ [\La_{ab}]*[H_{c}]=16[\La_{ab}]*[\omega_{c}]*[H_{c}]\\
& &=8[\La_{ac}]*[\La_{cb}]*[H_{c}]=128[\La_{ac}]*[\La_{cb}]*[\omega_{a}]*[H_{a}]\\
& &=64[\La_{ac}]*[\La_{ca}]*[\La_{ab}]*[H_{a}]=0.
\end{eqnarray*}
If $d=2$. Notice that $[\La_{ab}]=[S_{ab}(1,1)]$. By Remark 4.1.1
of \cite{DN3} and the fact that
$[\omega_{a}*S_{ab}(m,n)]=[S_{ab}(m,n)*\omega_{a}]$ in
$\bar{A}(V_{L}^{+})$ for $m,n\geq 1$, we have
\begin{equation}[S_{ab}(m+1,n)]+[S_{ab}(m,n)]=0.\label{ea30}
\end{equation} By the proof of Lemma 6.1.2 of \cite{DN3}, we know that
\begin{equation}
[H_{a}]=-9[S_{aa}(1,3)]-\frac{17}{2}[S_{aa}(1,2)]+\frac{1}{2}[S_{aa}(1,1)].\label{ea31}
\end{equation}
Direct calculation yields
$$
[S_{ab}(1,1)]*[S_{aa}(1,3)]=h_{b}(-1)h_{a}(-3)h_{a}(-1)^{2},$$
$$
[S_{ab}(1,1)]*[S_{aa}(1,2)]=h_{b}(-1)h_{a}(-2)h_{a}(-1)^{2},$$
$$
[S_{ab}(1,1)]*[S_{aa}(1,1)]=h_{b}(-1)h_{a}(-1)h_{a}(-1)^{2}.$$
Here we have used (\ref{ea30}). Then (\ref{ea5.25}) immediately
follows from Lemma 4.2.1 of \cite{DN3}, (\ref{ea30}) and
(\ref{ea31}).
 The proof is
complete. \qed

For $0\neq \al\in L$, let $\{h_{1},\cdots,h_{d}\}$ be  an
orthonormal basis of ${\f h}$ such that $h_{1}\in\C\al$. Define
$$
[\bar{B}_{\al}]=2^{|\al|^{2}-1}\frac{128}{9}[H_{1}]*[E^{\al}]. $$
We also set $[\bar{B}_{0}]=\frac{128}{9}[H_{1}]$.

\begin{lem}\label{l5.9}
 The subalgebra $A_{H}$ of $\bar{A}(V_{L}^{+})$ spanned by
$[\bar{B}_{\al}]$, $\al\in L$ is an ideal of $\bar{A}(V_{L}^{+})$
isomorphic to ${\C}[\widehat{L}/K]/J$.
\end{lem}

Let
$$
\widehat{A}(V_{L}^{+})=\bar{A}(V_{L}^{+})/A_{H}.$$

\begin{lem}\label{l5.10} Any $\widehat{A}(V_{L}^{+})$-module is
completely reducible. That is, $\widehat{A}(V_{L}^{+})$ is a
semisimple associative algebra.
\end{lem}

\pf Let $M$ be an $\widehat{A}(V_{L}^{+})$-module. For $\al\in L$,
by \cite{DN2} $M$ is a direct sum of irreducible
$A(V_{\Z\al}^{+})$-modules. Following the proof of Lemma 6.1 of
\cite{AD} one can prove that the image of any vector from
$M(1)^{+}$  in $\widehat{A}(V_{L}^{+})$ is semisimple on $M$. By
Table 1 of \cite{AD}, we can assume that
$$
M=\bigoplus_{\la\in{\f h/(\pm 1)}}M_{\la},$$ where $M_{\la}=\{w\in
M| \ [\frac{1}{2}h(-1)^{2}1]w=\frac{1}{2}(\la,h)^{2}w, h\in{\f
h}\}$. So $\omega_{a}w=\frac{1}{2}(\la,h_{a})^{2}w$, for $w\in
M_{\la}$. By (\ref{equation6}) and (\ref{equation7}), we have
$$\La_{ab}w=(\la,h_{a})(\la,h_{b})w, $$
for $a\neq b$, $w\in M_{\la}.$ For any $u\in M_{\la}$, $\la\neq
0$, set $M(u)=\sum\limits_{\al\in L}\C[E^{\al}]u$. By
(\ref{equation8})-(\ref{equation9}) and
(\ref{equation11})-(\ref{equation13}), if $[E^{\al}]u\neq 0$, then
$\al\in\De(\la)$ or $-\al\in \De(\la)$, where $\De(\la)=\{\al\in
L| \ |\la+\al|^{2}=|\la|^{2}\}$. So
$$M(u)=\bigoplus_{\al\in\De(\la)}\C[E^{\al}]u.$$
Since $L$ is positive-definite, there are finitely many $\al\in L$
which belong to $\De(\la)$. Thus $M(u)$ is finite-dimensional.
Similar to the proof of  Lemma 6.4 of \cite{AD}, we can deduce
that $\La_{ab}M(u)\subseteq M(u)$, $\omega_{a}M(u)\subseteq M(u)$.
By Proposition 5.4 of \cite{AD},
$[E^{\al}]*[E^{\be}]=[x]*[E^{\al+\be}]$ for some $x\in M(1)^{+}$.
We deduce that $M(u)$ is an $\widehat{A}(V_{L}^{+})$-submodule of
$M$. Suppose $[E^{\al}]u\neq 0$, for some $\al\in\De(\la)$. If
$(\al,\al)=2$, then by (\ref{equation10}), we have $0\neq
[E^{\al}][E^{\al}]u\in{\C}u$. If $(\al,\al)=2k\neq 2$. Let
$\{h_{1},\cdots,h_{d}\}$ be an orthonormal basis of ${\f h}$ such
that $h_{1}\in\C\al$. By the fact that
$[H_{1}]=[J_{1}]+[\omega_{1}]-4[\omega_{1}^{2}]=0$ and
(\ref{equation9}) we know that $[\omega_{1}]u=\frac{k}{4}u.$ Then
by Lemma 5.5 of \cite{DN2}, we have
$$
[E^{\al}][E^{\al}]u=\frac{2k^{2}}{(2k)!}(k^{2}-1)(k^{2}-2^{2})\cdots(k^{2}-(k-1)^{2})u\neq
0.$$
  Therefore $M(u)$ is
irreducible. We prove that $M$ is a direct sum of
finite-dimensional irreducible $\widehat{A}(V_{L}^{+})$-module.
\qed

\begin{theorem}\label{t5.11}
$V_{L}^{+}$ is a rational vertex operator algebra.
\end{theorem}

\pf By Lemmas \ref{l5.4}-\ref{l5.5}, \ref{l5.7}, \ref{l5.9} and
\ref{l5.10}, we know that $A(V_{L}^{+})$ is semisimple as $\C[\hat L/K]/J$
and $A^t\otimes \C[\hat L/K]/J$ are  semisimple.  Then the
theorem follows from Theorem \ref{t2.11}. \qed


\begin{thebibliography}{ABCDE}
\bibitem[A1]{A1} T. Abe, Rationality of the vertex operator algebra
$V_L^+$ for a positive definite even lattice $L,$ {\em Math. Z.} {\bf
249} (2005), 455-484.
\bibitem[A2]{A2} T. Abe, A $Z_2$-orbifold model of
the symplectic fermionic vertex operator superalgebra,
math.QA/0503472.
\bibitem[ABD]{ABD} T. Abe, G. Buhl and C. Dong, Rationality, regularity and
$C_2$-cofiniteness, {\em Trans. AMS.} {\bf 356} (2004), 3391-3402.
\bibitem[AD]{AD} T. Abe and C. Dong,
Classification of irreducible modules for the vertex operator algebra $V_L^+L:$ general case, {\em J. Algebra} {\bf  273} (2004), 657-685.
\bibitem[Ad]{Ad} D. Adamovi\'c, Classification of irreducible modules of
certain subalgebras of free boson vertex algebra, {\em  J.
Algebra} {\bf 270} (2003), 115-132.
\bibitem[B]{B}
R. Borcherds, Vertex algebras, Kac-Moody algebras, and the Monster,
{\it Proc. Natl. Acad. Sci. USA} {\bf 83} (1986), 3068-3071.
\bibitem[Bu]{Bu} G. Buhl, A spanning set for VOA modules,
{\em  J. Algebra} {\bf  254} (2002), 125-151.
\bibitem[D1]{D1} C. Dong, Vertex algebras associated with
even lattices, {\em J. Algebra} {\bf 160} (1993), 245-265.
\bibitem[D2]{D2} C. Dong, Representations of the moonshine module
vertex operator algebra, {\em Contemp. Math.} {\bf 175} (1994),
27-36.
\bibitem[DGH]{DGH} C. Dong, R. Griess Jr. and G. Hoehn,
Framed vertex operator algebras, codes and the moonshine module,
{\em Commu. Math. Phys.} {\bf 193} (1998), 407-448.
\bibitem[DJ1]{DJ1} C. Dong and C. Jiang, Bimodules associated to
vertex operator algebras, math.QA/0601626.
\bibitem[DJ2]{DJ2} C. Dong and C. Jiang, Representation theory of
vertex operator algebras  {\em Contemp. Math.,} to appear,
math.QA/0603588.
\bibitem[DLTYY]{DLTYY} C. Dong,  C.  Lam, K. Tanabe, H.  Yamada and K. Yokoyama, $\Bbb Z_3$ symmetry and $W_3$ algebra in lattice vertex operator algebras,
{\em  Pacific J. Math.}  {\bf 215} (2004), 245-296.
\bibitem[DL]{DL}  C. Dong and J. Lepowsky, Generalized Vertex
Algebras and Relative Vertex Operators, {\em Progress in Math.} Vol. 112,
Birkh\"{a}user, Boston 1993.
\bibitem[DLM1]{DLM1} C. Dong, H. Li and G. Mason, Regularity of rational vertex operator algebras, {\em Advances. in Math.}
{\bf 132} (1997), 148-166.
\bibitem[DLM2]{DLM2} C. Dong, H. Li and G. Mason,  Twisted representations of
vertex operator algebras, {\em Math. Ann.} {\bf 310} (1998),
571-600.
\bibitem[DLM3]{DLM3} C. Dong, H. Li and G. Mason, Vertex operator algebras and associative algebras, {\em J. Algebra}
{\bf 206} (1998), 67-96.
\bibitem[DLM4]{DLM5} C. Dong, H. Li and G. Mason,
Modular invariance of trace functions in orbifold theory and generalized
moonshine, {\em Commu. Math. Phys.} {\bf 214} (2000), 1-56.
\bibitem[DM1]{DM1}C. Dong and G. Mason, On quantum Galois theory,
{\em Duke Math. J.}, {\bf 86} (1997), 305-321.
\bibitem[DM2]{DM2} C. Dong and G. Mason,  Rational vertex operator algebras and the
effective central charge, {\em Internat. Math. Res. Notices} {\bf 56} (2004),  2989-3008.
\bibitem[DM3]{DM3} C. Dong and G. Mason, Integrability of $C_2$-cofinite vertex operator algebras, {\em Internat. Math. Res. Notices} {\bf 2006}
(2006), Art. ID 80468, 15 pp.
\bibitem[DMZ]{DMZ} C. Dong, G. Mason and Y. Zhu, Discrete series of the
Virasoro algebra and the moonshine module, {\em Proc. Symp. Pure. Math. American Math. Soc.} {\bf 56} II (1994), 295-316.
\bibitem[DN1]{DN1} C. Dong and K. Nagatomo, Classification of
irreducible modules for the vertex operator algebra $M(1)^+$,
{\em  J. Algebra} {\bf 216} (1999), 384-404.
\bibitem[DN2]{DN2} C. Dong and K. Nagatomo, Representations of vertex
operator algebra $V_L^+$ for rank one lattice
$L,$ {\em Comm. Math. Phys.} {\bf  202} (1999), 169-195.
\bibitem[DN3]{DN3} C. Dong and K. Nagatomo,
 Classification of irreducible modules for the vertex operator algebra $M(1)^+$. II. Higher rank, {\em J. Algebra} {\bf 240} (2001), 289-325.
\bibitem[DZ]{DZ} C. Dong and Z. Zhao,  Twisted representations of vertex
operator superalgebras, {\em Commu. Contemp. Math.} {\bf 8} (2006),
101-122.
\bibitem[FB]{FB} E. Frenkel and D. Ben-Ziv, Vertex Algebras and Algebraic
Curves, {\em Mathematical Surveys and Monographs} {\bf 88}, AMS. 2001.
\bibitem[FHL]{FHL} I. B. Frenkel, Y. Huang and J. Lepowsky, On
axiomatic approaches to vertex operator algebras and modules,
{\it Memoirs American Math. Soc.} {\bf 104}, 1993.
\bibitem[FLM]{FLM} I. B. Frenkel, J. Lepowsky and A. Meurman,
Vertex Operator Algebras and the Monster, {\em Pure and Applied
Math.} Vol. {\bf 134}, Academic Press, 1988.
\bibitem[FZ]{FZ} I. Frenkel and Y. Zhu, Vertex operator algebras associated to
representations of affine and Virasoro algebras, {\em  Duke Math. J}
{\bf  66} (1992), 123-168.
\bibitem[GN]{GN} M. Gaberdiel, and A. Neitzke, Rationality, quasirationality,
and finite W-algebras, {\em Comm. Math. Phys.} {\bf 238} (2003), 305-331.
\bibitem[H]{H} Y. Huang, Vertex operator algebras, the Verlinde conjecture, and modular tensor categories, {\em  Proc. Natl. Acad. Sci. USA} {\bf 102} (2005), 5352-5356.
\bibitem[KR]{KR} V. Kac and A. Raina, Highest Weight Representations of Infinite
Dimensional Lie Algebras, World Scientific, Adv. Ser. In Math. Phys.,
Singapore, 1987.
\bibitem[KW]{KW} V. Kac and W. Wang,
Vertex operator superalgebras and their representations, {\em
Contemp. Math. Amer. Math. Soc.} {\bf 175} (1994), 161-191.
\bibitem[KL]{KL} M. Karel and H. Li, Certain generating subspace for vertex
operator algebras, {\em J. Algebra} {\bf 217} (1999), 393-421.
\bibitem[KMY]{KMY} M. Kitazume, M.  Miyamoto and H. Yamada,
Ternary codes and vertex operator algebras, {\em  J. Algebra} {\bf
223} (2000), 379-395.
\bibitem[LL]{LL} J. Lepowsky and H. Li,
Introduction to Vertex Operator Algebras and Their Representations,
{\em  Progress in Mathematics} Vol. {\bf  227}, Birkhäuser Boston, Inc., Boston, MA, 2004.
\bibitem[L]{L} H. Li, Some finiteness properties of regular vertex
operator algebras, {\em J. Algebra} {\bf 212} (1999), 495-514.
\bibitem[M]{M} M. Miyamoto, Representation theory of code vertex operator
algebra, {\em  J. Algebra} {\bf  201} (1998), 115-150.
\bibitem[NT]{NT} K. Nagatomo and A. Tsuchiya, Conformal field theories associated to regular chiral vertex operator algebras I, Theories over the projective
line, {\em  Duke Math. J.} {\bf 128} (2005), 393-471.
\bibitem[TY]{TY} K. Tanabe and H. Yamada,
The fixed point subalgebra of a lattice vertex operator algebra by an automorphism of order three, math.QA/0508175.
\bibitem[W1]{W1} W. Wang, Rationality of Virasoro vertex operator algebras,
{\em Internat. Math. Res. Notices} {\bf  7} (1993), 197-211.
\bibitem[W2]{W2} W. Wang, Classification of irreducible modules of $W_3$ algebra with $c=-2$, {\em Comm. Math. Phys.}  {\bf 195} (1998), 113-128.
\bibitem[Z]{Z} Y. Zhu, Modular invariance of characters of vertex operator algebras,
{\em J. Amer, Math. Soc.} {\bf 9} (1996), 237-302.
\end{thebibliography}
\end{document}